%
%
\input amstex.tex
\documentstyle{amsppt}
\baselineskip=18truept
\voffset=0.3truein      
\hsize=6.0truein        
\vsize=8.5truein        
\parskip=5truept        
\overfullrule=0pt       
\TagsOnRight
\refstyle{A}
\widestnumber\key{He2}
\document
\font \myfont =cmbx10 scaled 1200
\font \nilsfont =cmbx10 scaled 1000
\input xypic
\define \s {\sigma}
\redefine \k {\kappa}
\redefine \l {\lambda}
\redefine \L {\Lambda}
\define \e {\epsilon}
\redefine \v {\varepsilon}
\define \p {\phi}
\redefine \P {\Phi}
\define \r {\rho}
\define \C {\Bbb C}
\define \N {\Bbb N}
\define \Z {\Bbb Z}
\define \R {\Bbb R}
\define \Hom {\text {Hom}}
\define \supp {\text {supp}}
\redefine \c {\cdot}
\redefine \H {\Cal H}
\define \F {\Cal F}
\define \a {\goth a}
\topmatter
\leftheadtext{Nils Byrial Andersen}
\rightheadtext{Fourier transform of Schwartz functions}
\pageno=1
\title On the Fourier transform of Schwartz functions on Riemannian Symmetric Spaces.\endtitle
\author by Nils Byrial Andersen.
\endauthor
\address
Department of Mathematics,
Aarhus University,
Ny Munkegade 118,
DK-8000 Aarhus C,
Denmark
\endaddress
\email
byrial\@imf.au.dk
\endemail
\abstract
Consider the (Helgason-) Fourier transform on a Riemannian symmetric space $G/K$.
We give a simple proof of the $L^p$-Schwartz space isomorphism theorem ($0 <p \le 2$) for $K$-finite
functions. The proof is a generalization of J.-Ph.\,Anker's proof for $K$-invariant functions.
\endabstract
\endtopmatter

\flushpar
{\myfont 1. Introduction.}

Let $G/K$ be a Riemannian symmetric space, where $G$ is a connected, non-compact
semisimple Lie group with finite center, and $K$ is the maximal compact subgroup fixed by a
Cartan involution. Let $\goth g = \goth k \oplus \goth p$ be the Cartan decomposition of $\goth g$,
the Lie algebra of $G$,
and let $\a$ be a maximal abelian subspace of $\goth p$. Let $M = Z _K (\goth a)$, then $K/M$ is a symmetric space.

Let as usual $\rho$ denote half the sum of the positive roots, and let $W$ denote the Weyl group. Let $\v\ge 0$.
Let $\goth a ^*$ denote the dual of $\goth a$, let  $C^{\v \r}$ be the convex hull of the
set $W \c \v\r$ in $\a ^*$,
and let $\goth a ^* _\v =\goth a ^* +iC^{\v \rho}$ be the tube with basis $C^{\v \rho}$ in the complex
dual $\goth a ^* _\C$.

Let $\H$ denote the (Helgason-) Fourier transform on $G/K$, and consider the
$L^p$-Schwartz spaces $\Cal S ^p (G/K),\,0 <p \le 2$, and the (semi-classical) Schwartz spaces $\Cal S (\goth a ^* _\v \times K/M)$.
The Schwartz space isomorphism theorem ([Eg,\,Theorem\,4.1.1]) states that:

\proclaim {Theorem 1}
Let $0 <p \le 2$ and $\v = \frac 2p - 1$.
The Fourier transform $\H$ is a topological isomorphism between $\Cal S ^p (G/K)$ and
$\Cal S (\goth a ^* _\v \times K/M)$, with the usual inverse.
\endproclaim

J.-Ph.\,Anker gave in [An] a simple and beautiful proof of Theorem 1 when restricted to $K$-invariant functions,
in which case the Fourier transform reduces to the spherical transform.
In Section 2 of this note we extend Anker's proof to $K$-finite functions for the real hyperbolic spaces (the rank 1 case), and in Section 3
we sketch how to prove the general case.

{\bf Remark:}
This note consists of two slightly reworked Chapters from
my "Progress Report": "On the Fourier transform on real Hyperbolic Spaces",
Aarhus University, 1995. Recently, two articles on the same subject were published:

J. Jana and P. Sarkar, {\it{On the Schwartz space isomorphism theorem for rank one symmetric space}},
Proc. Indian Acad. Sci. Math. Sci., {\bf 117} (2007), no. 3, 333--348, and

J. Jana, {\it{On the Schwartz space isomorphism theorem for the Riemannian symmetric spaces}}, arXiv:1002.4855.

Further material may also be found in Jana's Thesis: "Isomorphism of Schwartz
spaces under Fourier transform", Indian Statistical Institute, Kolkata,
July 2008.

The proofs by Jana and Sarkar use a reduction to the $K$-invariant result by Anker;
indeed they show that the various spaces of a fixed $K$-type are isomorphic to the similar space of trivial $K$-type.
The proof in the present note is a more straightforward generalization of Anker's proof,
with the use of a (generalized) Abel (or Radon) transform and a cut-off function
to show that the inverse transform is continuous from the Paley--Wiener space to the space of smooth functions with compact support on $G/K$.
Unfortunately, I also conclude that the proof cannot be generalized to the general case without restriction to
$K$-finite functions.
However, I hope, and think, that the present notes could form a nice supplement to the papers and Thesis mentioned above.

{\it Acknowledgments:}
As mentioned above, these notes were originally written during my Ph.D.-studies at Aarhus
University, Denmark. It is a pleasure to give belated thanks to my supervisors Henrik Stetk\ae r
and Henrik Schlichtkrull. Also it is a pleasure to thank Jean-Philippe Anker for useful discussions annd suggestions
concerning the proof during a short visit to Nancy in 1995.

\flushpar
{\myfont 2. The rank one case.}

In the following $c$ will denote (possibly different) positive constants.

\flushpar
{\nilsfont 2.1. Notation and preliminaries.}

Let $G$ be a connected semisimple Lie group with finite center, and let $\theta$ be a Cartan involution of $G$.
Then the fixed point group $K= G^\theta$ is a maximal compact subgroup. Let $\goth g$ and
$\goth k$ denote the respective Lie algebras, we then have a Cartan decomposition of $\goth g$ given by:
$\goth g = \goth k \oplus \goth p$. The Killing form on $\goth g$ induces an $Ad K$-invariant scalar
product on $\goth p$, and hence a $G$-invariant Riemannian metric on $X=G/K$.
With this structure, $X=G/K$ becomes a Riemannian globally symmetric space of the noncompact type.

Fix a maximal abelian subspace $\goth a$ of $\goth p$, and
let $\goth g = \goth k \oplus \goth a \oplus \goth n$ be an Iwasawa
decomposition of $\goth g$. This induces the Iwasawa decomposition $G = KAN$, where $A$ and $N$ are the Lie groups
corresponding to $\goth a $ and $\goth n$ respectively.
In this Section we will assume $G$ to be of rank one, i.e.,\,$\dim \goth a = 1$.
In the next Section we will sketch how to remove this condition.
We identify $\goth a$ and it's dual $\goth a ^*$ with
$\R$ (for the unique positive root, $\alpha \in \goth a ^*$,
we make the identification $\alpha =1$, for the unique element
$ H$ in $\goth a$ such that $\alpha (H) =1$, we put $H=1$. Then we define $a_t = \exp t H$).
Thus every $g \in G$ can be written as $k a_tn$, where $k\in K,t \in \R$ and $n \in N$ are unique.
We will denote the Iwasawa projections on the
$K$-part and $A$-part by $\kappa (g) = \kappa (ka_tn) = k$ and $H (g)= H(ka_tn) =t$.
Furthermore we will consider the "reverse" Iwasawa
decomposition, namely: $G = NAK$, where we will denote the projection onto the $A$-part by:
$A (g) =A(na_tk) =t$.
Remark that
$a_{-t} = a_t ^{-1}$ and $H(g) = - A (g ^{-1})$.

Define $A_+ = \{ a_t \in A | t >0\}$ and $\overline {A _+} = \{ a_t \in A | t \ge 0\}$,
corresponding to the open and closed
positive Weyl chambers, then we have the Cartan (or Polar) decomposition
of $G$ given by: $G = K \overline {A _+}K$, that is, every element
$g \in G$ can be written as $k_1 a_t k_2$, where $t \in \overline {\R_+}$
is unique and $k_1,k_2 \in K$.
We will define $|g| = |k_1 a_t k_2| = t$.
For the rank one case we have the basic estimate:
$|H(g)| \le  |g|$
(Given a finite-dimensional irreducible representation of $G$,
$H(g)$ and $|g|$ can be found using a normalized
highest weight vector). In the Cartan decomposition
the Haar measure is given by:
$$
\int _G f(x) dx = c\int _K \int _{\R _+} \int _K f (k_1 a_t k_2) \sinh ^n t  dk_1 dt dk_2,
$$where $n = \dim \goth n$.
In the hyperbolic case, $G = SO (p,1),\, K= SO(p)$, we get: $n = p-1$.
In $\goth a ^* \cong \R$ we define the element
$\r (H) = \frac 12 tr (ad H_{|\goth n }),\, H \in \goth
a$, or under the identification with $\R$: $\r  = \frac 12 n$.
We have the obvious estimate: $ 0 \le \sinh ^n t
\le c e ^{2\r t}$.

Define $M = Z _K (\goth a)$, then $B=K/M$ is a symmetric space.
We will sometimes identify functions on $X=G/K$ ($B = K/M$) as
right-$K$-invariant (right-$M$-invariant) functions on $G$ ($K$).
Then the invariant measure on $X$ is given by:
$$
\int _X f(x) dx = c\int _K \int _{\R _+}  f (k a_t K) \sinh ^n t  dk dt ,
$$
for some constant.
This is the Cartan decomposition of the measure on $X$.
Define $A(x,b)= A(gK,kM) =A(k^{-1} g)$.
\proclaim {Definition 2.1.1 (The Fourier transform on $G/K$)}
For $f\in C^\infty _c ( G/K)$, we define the Fourier transform by:
$$
\Cal H f(\nu,b) = \hat f (\nu,b) = \int _X f(x) e^{(-\nu + \rho ) (A(x,b))} dx,
$$for all $\nu \in \Bbb C , \, b \in B$.
\endproclaim

{\bf Remarks:}
Let $f \in C^\infty _c (K\backslash G/K)$. Then: ($k\in K$)
$$\split
\int _X f(x) e^{(-\nu + \rho ) (A(x,b))} dx &= \int _X f(k \c x) e^{(-\nu + \rho ) (A(x,b))} dx\\
&= \int _X f( x) e^{(-\nu + \rho ) (A(x,k \c b))} dx.
\endsplit\
$$
Since $K$ acts transitively on $B$,
we see that $\hat f$ is independent of $b\in B$.
Integrating over $K$ we get: ($b =kM$)
$$
=\int _K \int _X f(x) e^{(-\nu + \rho ) (A(x,kM))}dx dk = \int _X f(x) \int _K e^{(-\nu + \rho ) (A(x,kM))}dk dx.
$$
Let $x = gK$, we then recognize the spherical function $\varphi _{-\nu},\, \nu \in \C$ on $G$:
$$
\varphi _{-\nu} (g)=  \int _K e^{(-\nu + \rho ) (A(kg)}dk,
$$
see [He1,\,Theorem\,4.3].
For $f \in C^\infty _c (K\backslash G/K)$,
the Fourier transform thus reduces to the spherical transform:
$$
\Cal H f(\nu) = \int _X f(x)\varphi _{-\nu} (x) dx.
$$
treated in [An], [He1], [GV], etc.

From [He2,\,Chapter\,III,\,\S1,\,\S5],
we get Theorems 2.1.2, 2.1.3 and 2.1.5, where $c(\c)$ is the Harish-Chandra $c$-function:
\proclaim {Theorem 2.1.2 (The inversion formula)}
Let $f \in C^\infty _c ( X), \, x \in X$. Then:
$$
f(x) =c\int _{\Bbb R _+ \times B} e^{(i\nu + \rho ) (A(x,b))} \hat f (i\nu,b) {|c(i\nu)|}^{-2} d\nu db.
$$
\endproclaim
\proclaim {Theorem 2.1.3 (The Plancherel formula)}
Let $f_1, \, f_2 \in C^\infty _c ( X)$, then:
$$
\int _X f_1 (x) \overline {f_2 (x)} dx =
c\int _{\Bbb R _+ \times B} \hat f_1 (i\nu,b) \overline { \hat f_2 (i\nu,b)} {|c(i\nu)|}^{-2} d\nu db.
$$
The Fourier transform extends to an isometry of $L^2 (X)$ onto $L^2 (i\Bbb R _+ \times B, c{|c(i\nu)|}^{-2} $
$d\nu db)$.
\endproclaim
A $C^\infty$-function $\psi (z,b)$ on $\Bbb C \times B$, holomorphic in $z$,
is called a holomorphic function of uniform exponential type
$R$, if there exists a constant $R \ge 0$ such that for each $N \in \Bbb N$ we have:
$$
\sup _{z \in \Bbb C , b\in B } e^{-R | Re z | } {(1 + |z| )} ^N | \psi (z,b) | < \infty.
$$
\proclaim {Definition 2.1.4}
The space of holomorphic functions of uniform exponential type $R$ will be denoted $\Cal H ^R (\Bbb C \times B)$. Furthermore denote by
$\Cal H  (\Bbb C \times B)$ their union over all $R > 0$.
Let $\Cal H _e (\Bbb C \times B)$ denote the space of functions
$\psi \in \Cal H  (\Bbb C \times B)$
satisfying the symmetry condition (SC1):
$$
\int _B e^{(-\nu + \rho ) (A(x,b))} \psi (-\nu , b) db =
\int _B e^{(\nu + \rho ) (A(x,b))} \psi (\nu , b) db,\quad\nu \in \Bbb C ,\, x \in X.
$$
\endproclaim
We usually call $\Cal H  (\Bbb C \times B)$ the Paley-Wiener space.
For $\psi \in \Cal H _e (\Bbb C \times B)$ independent of $b$, (SC1) reduces to the
symmetry condition: $\varphi _{-\nu} (x) \psi (-\nu) =
\varphi _{\nu} (x) \psi (\nu)$, where $\varphi _\nu $ is the spherical function
indexed by $\nu$. Since $\varphi _{-\nu} = \varphi _\nu$ this
again reduces to $\psi$ being an even function of $\nu$, that is, the usual symmetry condition for the spherical transform (Weyl group
invariance), see [He1, p. 450].

Consider $f \in C_c ^\infty (X)$.
By the Cartan (polar) decomposition, the polar distance from the point $x = gK = ka_t$ to the origin $x_o = eK$
is $|x| = |t|$.
Let $R>0$. We say that $\supp (f) \subset \bar B (0,R)$
if and only if the function $f$ has support inside $K \times  \bar B (0,R)$ or
if and only if $f(x) = 0$ for $|x| > R$.
\proclaim {The Paley-Wiener Theorem 2.1.5}
The Fourier transform is a bijection of the space $C_c ^\infty (X)$
onto the space $\Cal H _e (\Bbb C \times B)$, the inverse
transform being given by Theorem 2.1.2. Moreover
$\supp (f) \subset \bar B (0,R)$ if and only if $\hat f \in \Cal H _e  ^R (\Bbb C \times B)$.
\endproclaim
We decompose the Fourier transform into two transforms.
Let $f \in C_c ^\infty (G/K)$, then the Radon transform is defined as follows:
$$
\Cal R f (t,k) = e^{\rho t} \int _N f(ka_t n) dn.
$$
\proclaim {Proposition 2.1.6}
The Radon transform $\Cal R$ maps $C_c ^\infty (G/K)$ into $C_c ^\infty (\Bbb R \times B)$.
If $\supp f \subset \bar B (0, R)$, then $\supp \Cal R f \subset
\bar B (0,R) \times B$.
\endproclaim
\demo {Proof}
Let $ k, k_1, k_2 \in K, n\in N , t\in \R$,
we then get: $|k a_{t} n| \ge |t|  = |k_1 a_{|t|}k_2|$, hence $\supp \Cal R f
\subset \supp f \times B$.
$\hfill \square$
\enddemo
{\bf Remark:} For $f \in C^\infty _c (K\backslash G/K)$,
the Radon transform reduces to the Abel transform, see [An].

Let $\p \in C_c ^\infty (\R \times B)$,
then the "classical" Fourier transform on $\Bbb R \times B$ is defined as:
$$
\Cal F \p (\nu ,b) = \int _{\Bbb R} \p (t,b) e^{-\nu t} dt ,\quad \nu \in i\Bbb R ,\, b \in B.
$$
Let $\psi$ be a nice function on $i\R \times B$,
then the "classical" inverse Fourier transform is defined by:
$$
{\Cal F}^{-1} \psi (t,b) =\frac 1{2\pi} \int _{\Bbb R} \psi (i\nu ,b) e^{i\nu t} d\nu ,\, t \in \Bbb R ,\, b \in B.
$$We then have: ($b = kM$)
$$
\hat f (\nu , kM) = \int _{\Bbb R} e^{-\nu t} \left \{ e^{\rho t} \int _N f(ka_t n) dn \right \} dt,
$$
i.e.,\,we have the following commutative diagram:
$$
\diagram
& { {\Cal H}_e (\Bbb C \times B) }   \\
 { C_c^\infty ( G /K)} \urto^{\Cal H} \rrto_{\Cal R}& &
{ \Cal R C_c^\infty ( G /K) \subset C_c ^\infty (\Bbb R \times B)}  \ulto_{\Cal F} \\
\enddiagram
$$
The commutativity is an easy consequence of the definitions of the various transforms. Furthermore we get:
\proclaim {Proposition 2.1.7}
The Radon transform $\Cal R$ is an isomorphism between $C_c ^\infty ( G/K)$ and $\Cal R C_c ^\infty (G/K) = {\Cal F}^{-1}
\Cal H _e (\Bbb C \times B) \subset C_c ^\infty (\Bbb R \times B)$. Moreover $\supp f \subset \bar B (0, R)$ if and only if $\supp \Cal R f
\subset \bar B (0, R) \times B$.
\endproclaim
\demo {Proof}
The Paley-Wiener Theorems above and below.
$\hfill \square$
\enddemo
{\bf Remark:} We have indirectly introduced a symmetry condition for functions in
$C_c ^\infty (\Bbb R \times B)$, namely
that the Fourier transformed function should be in $\Cal H _e (\Bbb C \times B)$.
For functions of a specific $K$-type ($K$
acting on the $B$-variable), this
symmetry condition becomes somewhat easier to describe, see later.
For the trivial $K$-type, or for $B$-invariant functions, it
reduces to: $g (-t) = g(t), g \in C_c ^\infty (\Bbb R)$, or $g$ even,
that is, as in [An], where the Abel transform maps $C_c ^\infty (
K\backslash G/K)$ onto even functions in $C_c ^\infty (\Bbb R)$.
\proclaim {Theorem 2.1.8 (A Paley-Wiener Theorem on  $\Bbb R \times B$)}
Let $\phi \in C_c ^\infty ( \Bbb R \times B )$ have support in $\bar B (0,R)\times B$,
and let $f (\nu ,b) = \int _{\Bbb R} \phi (t,b) e^{-\nu t} dt
 ,\, \nu \in i\Bbb R ,\, b \in B$.
Then $f \in C ^\infty ( \Bbb C \times B )$, and $f(\cdot , b)$
is an entire function for fixed $b$, such that for all $N \in\Bbb N$, we have:
$$
\sup _{z \in \Bbb C , b\in B } e^{-R | Re z | } {(1 + |z| )} ^N | f (z,b) | < \infty.
$$
Conversely, let $f \in C ^\infty ( \Bbb C \times B )$ satisfy the above. Then there exists a function
$\phi \in C_c ^\infty ( \Bbb R \times B )$,
with support in $\bar B (0,R)\times B$, such that $f = \Cal F \phi$.
\endproclaim
\demo {Proof}
An easy generalization of [Ru, Theorem 7.22].
$\hfill \square$
\enddemo
We will in the following sections need some estimates on the spherical functions
$\varphi _\nu$ introduced before and on the Harish-Chandra $c$-function.
\proclaim {Lemma 2.1.9}
The spherical functions are all bi-$K$-invariant, $\varphi _{-\nu} =\varphi _\nu$, and:
\roster
\item"i)" For $ t\ge 0,\,\exists c >0$ such that: $ e ^{-\r t} < \varphi _o (a_t) < c (1+t) e ^{-\r t}$.
\item"ii)" Let $\varepsilon \ge 0,\, |Re \nu |
\le \varepsilon \rho$ and $ t \ge 0 $. Then: $|\varphi _\nu (a_t)| \le c (1+t) e^{(\varepsilon -1) \rho t}$.
\endroster
For the $c$-function we have:
\roster
\item"iii)" Let $Re \nu \ge 0$.
$\exists \gamma \ge 0 $ such that: ${|c(\nu)|}^{-1} \le \gamma(1+|\nu |)^{\frac {n}{2}}$.
\endroster
\endproclaim
\demo {Proof}
i) [He1, Chap IV, Ex.B1; GV, Sect.4.6].
\flushpar
ii) Assume $\nu \ge 0$: Then:
$$
{\varphi}_{\nu} (a_t) \le e^{\nu t} {\varphi}_{o} (a_t)
\le c e^{\varepsilon \rho t}  (1+t) e^{-\rho t} =  c (1+t) e^{(\varepsilon -1) \rho t}.
$$
\flushpar
iii) Properties of the $\Gamma$-function, see [He1,\,Chapter IV,\,Proposition 7.2].
$\hfill \square$
\enddemo

\flushpar
{\nilsfont 2.2. Schwartz spaces and dense subspaces.}

We have come to the definition of the Schwartz spaces.
Let $U(\goth g)$ denote the
universal enveloping algebra of $\goth g$. The elements of $U(\goth g)$ act on $C^{\infty} (G)$ as differential operators on both
sides.
We shall write $f(D;g;E)$ for the action of
$(D,E) \in U(\goth g) \times U(\goth g)$ on $f \in C^{\infty} (G)$ at $g \in G$, more explicitly we have:
$$\split
f(D;g;E) &=\left . {\left ( \frac {\partial}{\partial s_1} \dots \frac {\partial}{\partial s_d}\frac {\partial}{\partial t_1}\dots \frac {\partial}{\partial t_e} \right )} \right | _
{s_1 =\dots =s_d = t_1 = \dots = t_e =0}\\ & \times f((\exp s_1 X_1)\cdots (\exp s_d X_d)g(\exp t_1 Y_1)\cdots (\exp t_e Y_e)).\endsplit
$$
if $D=X_1 \cdots X_d, \, E = Y_1\cdots Y_e \, (X_1,...,X_d,Y_1,...,Y_e \in \goth g)$.
\proclaim {Definition 2.2.1}
Let $0<p \le 2$. The $L^p$-Schwartz space $\Cal S ^p (G /K)$ is the space of all functions $f \in C^{\infty} ( G /K)$, such that
$$
\sup _{g \in G} {(1+ |g|)}^N {\varphi _o (g)}^{- \frac 2p} |f(D;g;E)| < \infty,
$$
for any $D,E \in U(\goth g)$, and any nonnegative integer $N$.
Here $\varphi _o$ is the spherical function with $\nu = 0$.
\endproclaim
The topology of $\Cal S ^p (G /K)$ is defined by the seminorms:
$$
\sigma _{D,E,N} ^p  (f) = \sup _{g \in G} {(1+ |g|)}^N {\varphi _o (g)}^{- \frac 2p} |f(D;g;E)|.
$$
{\bf Remarks} 1) Consider the Cartan decomposition of the integral of $X$, then:
$$
\int _X f(x) dx = \int _{K\times \R _+} f (ka_t) \sinh ^n t dk dt.
$$
From above we see that the factor ${\varphi _o (g)}^{- \frac 2p}$
will control the factor $\sinh ^n t$, and hence the definition
resembles the natural definition of a Schwartz space on $K \times \R _+$.
\flushpar
2) In fact we see that: $\Cal S ^p (G /K) \subset L^q (G /K)$ for $0 < p\le q \le 2$,
while $\Cal S ^p (G /K) \not\subset L^q (G /K)$ for
$0 < q <p \le 2$. For $0 < p\le q \le 2$ the inclusion above is continuous.
\flushpar
3) Obviously $C_c^\infty (G /K) \subset \Cal S ^p (G /K)$ for $0 <p\le 2$,
and the inclusion is continuous.
\proclaim {Lemma 2.2.2} Let $0 < p\le q \le 2$. Then:
\roster
\item "(i)" $\Cal S ^p (G /K)$ is a Fr\'echet space.
\item "(ii)" $C_c^\infty (G /K)$ is a dense subspace of $\Cal S ^p ( G /K)$.
\item "(iii)" $\Cal S ^p ( G /K)$ is a dense subspace of $\Cal S ^q ( G /K)$.
\item "(iv)" $\Cal S ^p ( G /K)$ is a dense subspace of $L^q ( G /K)$.
\item "(v)" $\Cal S ^p (K \backslash G /K)$ is a closed subspace of $\Cal S ^p (G /K)$ in the Fr\'echet topology.
\endroster\endproclaim
\demo {Proof}
See [GV, sect. 6.1, 7.8].
$\hfill\square$
\enddemo
For
$E \in U(\goth k)$ define: $E f(k) = f(k;E)$.
Let $\Omega _K$ denote the Casimir element of  $U(\goth k)$, then for the Laplace-Beltrami operator
$\Delta _B$ on $B$ we have (modolo a constant):
$\Delta _B f = f(\cdot ; \Omega _K) = f(\Omega _K ; \cdot)$ ($\Omega _K \in \goth Z (\goth k )$).
Furthermore $\left ( P \left (\frac {\partial}{\partial \nu} \right) , E \right )$ will denote
differentiation with $P \left (\frac {\partial}{\partial \nu} \right)$
on the $\nu$-variable and with $E$ on the $k$-variable.
\proclaim {Definition 2.2.3}
Fix $\varepsilon \ge 0$.
Let $i{\Bbb R}_\varepsilon = i{\Bbb R} +[-\varepsilon \rho ,\varepsilon \rho]$.
The Schwartz space $\Cal S (i{\Bbb R}_\varepsilon \times B)$
consists of all complex valued functions $f \in C^\infty (i{\Bbb R}_\varepsilon \times B)$
such that:
\roster
\item "(i)" For fixed $b \in B$, $f (\cdot , b)$ is holomorphic in the interiour of $i{\Bbb R}_\varepsilon \, (i{\Bbb R}^o _\varepsilon)$.
\item "(ii)" $f$ and all its derivatives extend continuously to $i{\Bbb R}_\varepsilon \times B$.
\item "(iii)" For any polynomial $P$, any $E \in U(\goth k)$ and nonnegative integer $N$ we have the estimate:
$$
\sup _{\nu \in i{\Bbb R}_\varepsilon , \, k \in K} {(1+ |\nu|)}^N
\left|{\left (P\left (\frac {\partial}{\partial \nu} \right) ,E \right )f(\nu ,k)}\right | < \infty.
$$
\endroster
\endproclaim
The topology of $\Cal S (i{\Bbb R}_\varepsilon \times B)$ is defined by the seminorms:
$$
\tau _{P,E,N} ^\varepsilon (f) = \sup _{\nu \in i{\Bbb R}_\varepsilon  , \, k \in K} {(1+ |\nu|)}^N
\left|{\left (P\left (\frac {\partial}{\partial \nu} \right) ,E  \right )f(\nu ,k)}\right | < \infty.
$$
{\bf Remarks:}
For $\varepsilon = 0$ condition i) is empty. For $\varepsilon > 0$ ii) and iii) are equivalent to:
\roster
\item "(iv)" For any polynomial $P$, any $E\in U(\goth k)$ and any nonnegative integer $N$ we have the estimate:
$$
\sup _{\nu \in i{\Bbb R}^o _\varepsilon , \, k \in K} {(1+ |\nu|)}^N
\left|{\left (P\left (\frac {\partial}{\partial \nu} \right) ,E  \right )f(\nu ,k)}\right | < \infty.
$$\endroster
To see this, observe that for fixed $k \in K$, ${\left (P\left (\frac {\partial}{\partial \nu} \right) ,E  \right )f(\nu ,k)}$ is bounded on
$i{\Bbb R}^o _\varepsilon$ for all polynomials
$P$. $i{\Bbb R}^o _\varepsilon$ is convex, hence by the Mean Value Theorem we get:
$$
| f(x_1,k) - f(x_2,k) | \le \{ \sup _{\nu \in i{\Bbb R}^o _\varepsilon}  |  \nabla f(\nu,k) | \}
| x_1- x_2 | ,\quad x_1 , x_2 \, \in i{\Bbb R}^o _\varepsilon ,
$$where $\nabla$ is the operator $\nabla = \frac {df}{d\nu}$,
which shows that $f(\c ,k)$ is uniformly continuous on $i{\Bbb R}^o _\varepsilon$, hence by density $f(\c,k)$ extends to a
uniformly continuous function on $i{\Bbb R} _\varepsilon$. This can also be done for all derivatives of $f$.
We also note that the space defined above is homeomorphic to the space of functions $f$ in
$C^\infty (i{\Bbb R}_\varepsilon \times B)$ satisfying i), ii) and:
\roster
\item "(v)" For any polynomial $P$ and any nonnegative integers $M$ and $N$ we have the estimate:
$$
\sup _{\nu , \, k \in K} {(1+ |\nu|)}^N
\left|{\left (P\left (\frac {\partial}{\partial \nu} \right) ,{\Omega _K}^M \right )f(\nu ,k)}\right | < \infty .
$$
\endroster
Here we define the topology by the seminorms:
$$
\tau _{P,M,N} ^\varepsilon (f) =\sup _{\nu , \, k \in K} {(1+ |\nu|)}^N
\left|{\left (P\left (\frac {\partial}{\partial \nu} \right) ,{\Omega _K}^M \right )f(\nu ,k)}\right | < \infty ,
$$
that is, restriction to powers of the Casimir element does not alter the space or topology,
see e.g. [Eg,\,p.193].

Let ${\Cal S (i{\Bbb R}_\varepsilon \times B)}_e$ denote the space of functions
$\psi \in {\Cal S (i{\Bbb R}_\varepsilon \times B)}$ satisfying
the symmetry condition (SC2):
$$
\int _B e^{(-\nu + \rho ) (A(x,b))} \psi (-\nu , b) db = \int _B e^{(\nu + \rho ) (A(x,b))} \psi (\nu , b) db,\quad
\nu \in i{\Bbb R}_\varepsilon ,\, x \in X.
$$
\proclaim {Lemma 2.2.4} Let $\varepsilon \ge 0$.
\roster
\item "(i)" ${\Cal S (i{\Bbb R}_\varepsilon \times B)}_e$ is a Fr\'echet space.
\item "(ii)" ${\Cal H} (\Bbb C \times B) $ is a dense subspace of ${\Cal S (i{\Bbb R}_\varepsilon \times B)}$.
\endroster\endproclaim
\demo {Proof}
(i) Clearly ${\Cal S (i{\Bbb R}_\varepsilon \times B)}_e$ is a closed subspace
of the Fr\'echet space ${\Cal S (i{\Bbb R}_\varepsilon \times B)}$,
hence it is a Fr\'echet space in the inherited topology.
\flushpar
(ii) Let $\Cal S (\Bbb R \times B)$ be the "classical" Schwartz space on $\Bbb R \times B$,
i.e.\,the space of functions  $f \in C^\infty (\Bbb R \times B)$ such that:
$$
\sup _{t \in \Bbb R ,\,  k \in K} \left | {(1 +|t|)}^N \left (
P\left (\frac {\partial}{\partial t} \right ) ,E \right ) f(t,k) \right | < \infty.
$$
for any polynomial $P$, $E \in U(\goth k )$ and $N \in \Bbb N$.
We consider the function $t \mapsto \cosh (\varepsilon \rho t)$.
Let $\Cal S _{\varepsilon \rho } (\Bbb R \times B)$ be the space of functions $f \in
C^\infty (\Bbb R \times B)$ such that:
$$
\sup _{t \in \Bbb R ,\,  k \in K} \left | {(1 +|t|)}^N \cosh (\varepsilon \rho t) \left (
P\left (\frac {\partial}{\partial t} \right ) ,E \right ) f(t,k) \right | < \infty.
$$
for any polynomial $P$, $E \in U(\goth k )$ and $N \in \Bbb N$.
Then $\Cal S _{\varepsilon \rho } (\Bbb R \times B)$ is a Fr\'echet space for the obvious topology. We now have two topological isomorphisms:
\flushpar
1) the map $f \mapsto \cosh (\varepsilon \rho \cdot ) f$ between $\Cal S _{\varepsilon \rho } (\Bbb R \times B)$ and the classical Schwartz space $\Cal S (\Bbb R \times B)$. (obvious)
\flushpar
2) the "classical" Fourier transform $\Cal F$ between $\Cal S _{\varepsilon \rho } (\Bbb R \times B)$ and $\Cal S (i{\Bbb R}_\varepsilon \times B)$.

Consider:
$$
\hat f (\nu ,b) = \int _{\Bbb R} f(t,b) e^{-\nu t} dt.
$$
The factor $e^{-\nu t}$ is bounded by $2\cosh (\varepsilon \rho t)$ for $ | Re \nu | \le \varepsilon \rho $, hence we can write the above as:
$$
\hat f (\nu ,b) = \int _{\Bbb R} 2\cosh (\varepsilon \rho t )
f(t,b) \frac 1{2\cosh (\varepsilon \rho t)} e^{-\nu t} dt .
$$
For fixed $k \in K$, we see by Morera's Theorem that $\hat f (\cdot , k)$
is holomorphic in $i{\Bbb R}^o _\varepsilon $. Now consider:
$$\split
&\left | {(1+ \nu)}^N \left ( P\left (\frac {\partial}{\partial \nu} \right ) ,E \right ) \hat f (\nu ,k)\right | \\
&= \left |\int _{\Bbb R} \left \{  {\left ( 1 + \left ( \frac {d}{dt} \right ) \right )}^N P(-t) Ef(t,k) \right \} e^{-\nu t} dt \right | \\
&\le c\sup _{t \in \Bbb R ,\,b \in B} \left |\cosh (\varepsilon \rho t ) {(1 +|t|)}^2  {\left ( 1 + \left ( \frac {d}{dt} \right ) \right )}^N
P(-t)Ef(t,k) \right | \\
&\le c\sup _{t \in \Bbb R ,\,  b \in B} \left |\cosh (\varepsilon \rho t ) {(1 +|t|)}^M \left ( \tilde P\left (\frac {\partial}{\partial t} \right ) ,E \right )
f(t,k)\right |
< \infty \, ,\endsplit\
$$
for some polynomial $\tilde P$ and $M \in \Bbb N$,
i.e.,\,$\hat f \in \Cal S (i{\Bbb R}_\varepsilon \times B)$
and the Fourier transform is continuous as an operator from $\Cal S _{\varepsilon \rho } (\Bbb R \times B)$
to $\Cal S (i{\Bbb R}_\varepsilon \times B)$.
Now consider the inverse Fourier transform:
$$
\check f (t,k) =\frac 1{2\pi} \int _{\Bbb R} f (i\nu ,k) e^{i\nu t} d\nu.
$$
Let $P$ and $Q$ be polynomials, then by Cauchy's Theorem,
shifting integral from $i\nu$ to $i\nu \pm \varepsilon \rho$,
we get:
$$
\split
&P(t) \cosh (\varepsilon \rho t ) \left ( Q \left (\frac {\partial}{\partial t} \right ) ,E \right ) \check f (t,k)=\\
&\frac 1{4 \pi} \int _{\Bbb R} P \left (i\frac {\partial}{\partial \nu} \right )
\{ Q(i\nu - \varepsilon \rho) Ef (i\nu - \varepsilon \rho ,k) +Q(i\nu + \varepsilon \rho) Ef (i\nu + \varepsilon \rho ,k)\} e^{i\nu t} d\nu \, , \endsplit\
$$
i.e.,\,
$$\split
&\sup _{t\in \Bbb R} \left | P(t) \cosh (\varepsilon \rho t) \left (Q\left (\frac {\partial}{\partial t} \right ) ,E \right ) \check
f(t,k) \right | \\
&\le c \sup _{\nu \in i{\Bbb R}_\varepsilon} \left|{(1+ |\nu|)}^2 {P\left (\frac {-\partial}{\partial \nu} \right) Q (\nu ) Ef(\nu
,k)}\right | < \infty \, ,
\endsplit\
$$
and we see that the inverse Fourier transform also is continuous.

Since $C_c ^\infty ({\Bbb R} \times B)$ is a dense subspace of $\Cal S({\Bbb R} \times B)$, we
get from above that $C_c ^\infty ({\Bbb R} \times B)$ is a dense subspace of
$\Cal S _{\varepsilon \rho } (\Bbb R \times B)$, and hence from the Paley-Wiener Theorem on $\Bbb R \times B$ we
conclude ii).
$\hfill\square$
\enddemo
Let in the following $\varepsilon \ge0$.
We want to show that $\Cal H _e (\Bbb C \times B)$ is a dense subspace of
${\Cal S (i{\Bbb R}_\varepsilon \times B)}_e$, that is, we have to
look at functions in each space satisfying the symmetry conditions (SC1-2).
We will make this problem a little easier by looking at $K$-types
for the left-regular representation of $K$ ($l$) on the Fr\'echet space
${\Cal S (i{\Bbb R}_\varepsilon \times B)}$ (acting on the second variable).
It is straightforward to check that this representation is a
smooth Fr\'echet representation of $K$. Let $\hat K$ denote the set of equivalence
classes of finite dimensional unitary irreducible representations
$(\delta , V_\delta )$ of $K$. Define $V_\delta ^M \equiv
\{v\in V_\delta | \delta (m) v = v ,\, m \in M\} $ and
$\hat K _M \equiv \{\delta \in \hat K | V_\delta ^M \ne
0 \}$. For $\delta \in \hat K$, let $d (\delta)$ and
$\chi _\delta$ denote the dimension and
character of $\delta$. Let ${\Cal S (i{\Bbb R}_\varepsilon \times B)}_\delta$
be the closed subspace of ${\Cal S (i{\Bbb R}_\varepsilon \times B)}$
consisting of functions of $K$-type $\delta$.
The continuous projection of ${\Cal S (i{\Bbb R}_\varepsilon \times B)}$ onto
${\Cal S (i{\Bbb R}_\varepsilon \times B)}_\delta$ is given by:
$$
P_\delta f (\cdot , b ) = f_\delta (\cdot , b ) =
d (\delta) \int _K \chi _\delta (k^{-1}) f (\cdot , k^{-1} \cdot b) dk.
$$[He1, Chapter IV, Lemma 1.7].
The space of $K$-finite functions in ${\Cal S (i{\Bbb R}_\varepsilon \times B)}$ is
given as:
$$
{\Cal S (i{\Bbb R}_\varepsilon \times B)}_K \equiv \{ f \in {\Cal S (i{\Bbb R}_\varepsilon \times B)} | \dim \text {span} \, l(K) f < \infty \}.
$$
Then $f \in {\Cal S (i{\Bbb R}_\varepsilon \times B)}_K \Leftrightarrow f = \sum _{\delta \in \hat K _M} f_\delta$, where the sum is finite.
For $\delta \in \hat K$ we will also consider the contragredient representation
$\check {\delta} \in \hat K$.
$\check {\delta} (k)$ can be identified as the operator in $\Hom (V_\delta ', V_\delta ')$
defined by $\check {\delta} (k) =  {\delta (k^{-1})}^t$,
where $V_\delta '$ is the dual space of $V_\delta$ and $t$ denotes
transpose. We remark that $\delta \in \hat K _M \Leftrightarrow \check {\delta} \in \hat K _M$.
We have the following important result, see [He1, Chapter V, Theorem 3.1]:
\proclaim {Theorem 2.2.5}
${\Cal S (i{\Bbb R}_\varepsilon \times B)}_K$ is a dense subspace of ${\Cal S (i{\Bbb R}_\varepsilon \times B)}$. Furthermore we have the expansion
$f = \sum _{\delta \in \hat K _M} f_\delta  = \sum _{\delta \in \hat K _M} f_{\check {\delta}}$,
where the sums are absolutely convergent.
\endproclaim
We easily see that $\Cal H (\Bbb C \times B)$ is invariant under $P_\delta$,
hence we will consider the subspace
$\Cal H (\Bbb C \times B) _\delta = \Cal H (\Bbb C \times B) \cap {\Cal S (i{\Bbb R}_\varepsilon \times
B)}_\delta$.
Since $P_\delta$ is continuous we see that
$\Cal H (\Bbb C \times B) _\delta$ is a dense subspace of
${\Cal S (i{\Bbb R}_\varepsilon \times B)}_\delta$.
Denote by $\Cal H (\Bbb C \times B) _{\delta ,e}$, respectively by
${\Cal S (i{\Bbb R}_\varepsilon \times B)}_{\delta ,e}$,
the set of functions in $\Cal H (\Bbb C \times B) _\delta$ respectively in
${\Cal S (i{\Bbb R}_\varepsilon \times B)}_\delta$
satisfying (SC2) (a function in $\Cal H (\Bbb C \times B) _\delta$
that satisfies (SC2) will automatically by
holomorphicity satisfy (SC1)).
We want to show that $\Cal H (\Bbb C \times B) _{\delta ,e}$ is a dense subspace of
${\Cal S (i{\Bbb R}_\varepsilon \times B)}_{\delta ,e}$, for all $\delta \in  \hat K _M$,
and then by Theorem 2.2.5 conclude that
$\Cal H _e (\Bbb C \times B)$ is a dense subspace of ${\Cal S (i{\Bbb R}_\varepsilon \times B)}_e$.

Let $\delta \in \hat K _M$ act on $V_\delta$.
Looking at matrix entries, we can define the spaces $\Cal H (\Bbb C \times B ,\Hom (V_\delta ,V_\delta ))$
and $\Cal S (i\Bbb R _\varepsilon \times B,\Hom (V_\delta ,V_\delta ))$
of Paley-Wiener functions on $\Bbb C \times B$, respectively
Schwartz functions on $i\Bbb R _\varepsilon \times B$, taking values in $\Hom (V_\delta ,V_\delta
)$. We equip $\Cal S (i\Bbb R _\varepsilon \times B,\Hom (V_\delta ,V_\delta ))$ with the obvious topology.
Define $\Cal S (i\Bbb R _\varepsilon \times B) ^\delta \equiv \{F \in \Cal S (i\Bbb R _\varepsilon
\times B,\Hom (V_\delta ,V_\delta )) : F(\cdot , k\cdot b ) = \delta (k) F (\cdot , b ) \}$,
and define
$\Cal H (\Bbb C \times B ) ^\delta =\Cal S (i\Bbb R _\varepsilon \times B) ^\delta \cap \Cal H (\Bbb C \times B ,\Hom (V_\delta ,V_\delta
))$.
For $f \in {\Cal S (i{\Bbb R}_\varepsilon \times B )}$, we define:
$$
P^\delta f (\cdot , b ) = f^\delta (\cdot , b ) = d (\delta) \int _K \delta (k) f (\cdot , k^{-1} \cdot b) dk.
$$
We easily see that $ P^\delta$ takes $\Cal H (\Bbb C \times B )$
into $\Cal H (\Bbb C \times B ) ^\delta$, and
${\Cal S (i{\Bbb R}_\varepsilon \times B )}$ into $\Cal S (i\Bbb R _\varepsilon \times B) ^\delta$.
Denote by $\Cal H (\Bbb C \times B ) ^\delta _e$, respectively by
$\Cal S (i\Bbb R _\varepsilon \times B) ^\delta _e$, the
set of functions $F$ in $\Cal H (\Bbb C \times B ) ^\delta$
and $\Cal S (i\Bbb R _\varepsilon \times B) ^\delta$ satisfying:
$$
\int _B e^{(-\nu + \rho ) (A(x,b))} F (- \nu ,b) db =
\int _B e^{(\nu + \rho ) (A(x,b))}  F (\nu ,b) db,\quad \nu \in i{\Bbb R}_\varepsilon,\, x \in X.
$$
\proclaim {Proposition 2.2.6} Let $\delta \in \hat K _M$.
\roster
\item "(i)" The map:
$$
Q:\, F(\nu ,b) \to Tr (F(\nu ,b))
$$is a homeomorphism of $\Cal S (i\Bbb R _\varepsilon \times B) ^\delta$ onto
${\Cal S (i{\Bbb R}_\varepsilon \times B)}_{\check \delta}$, with inverse
$f \to P^\delta f = f^\delta$.
Furthermore $Q$ takes $\Cal H (\Bbb C \times B ) ^\delta$
onto $\Cal H (\Bbb C \times B ) _{\check\delta}$.
Here $Tr$ denote trace in $\Hom (V_\delta ,V_\delta )$.
\item "(ii)" The image of $\Cal H (\Bbb C \times B ) ^\delta _e$,
respectively of
$\Cal S (i\Bbb R _\varepsilon \times B) ^\delta_e$,
under $Q$ is $\Cal H (\Bbb C \times B) _{\check \delta ,e}$,
respectively ${\Cal S (i{\Bbb R}_\varepsilon \times B)}_{\check \delta,e}$.
\item "(iii)" The maps:
$$\gather
P_\delta:\,{\Cal S (i{\Bbb R}_\varepsilon \times B)} \to {\Cal S (i{\Bbb R}_\varepsilon \times B)}_{ \delta} \quad \text {and} \\
P^\delta:\,{\Cal S (i{\Bbb R}_\varepsilon \times B)} \to {\Cal S (i{\Bbb R}_\varepsilon \times B)}^{ \delta}\endgather
$$
are continuous open surjections, and the images are closed in
${\Cal S (i{\Bbb R}_\varepsilon \times B)}$ and
$\Cal S (i\Bbb R _\varepsilon \times B,\Hom (V_\delta ,V_\delta ))$ respectively.
\endroster\endproclaim
\demo {Proof}
(i),\,(iii) As in [He1,p.395f].
\flushpar
(ii) Consider $f \in {\Cal S (i{\Bbb R}_\varepsilon \times B)}_{\check \delta,e}$.
In (SC2) replace $x$ by $k^{-1} \cdot x$, use $A(k^{-1} \cdot x ,b) = A(x, k\cdot b)$,
make the substitution $b \mapsto k^{-1} \cdot b$, multiply by $\delta (k)$
and integrate over $K$. The result is:
$$
\int _B e^{(-\nu + \rho ) (A(x,b))} f^\delta (- \nu ,b) db = \int _B e^{(\nu + \rho ) (A(x,b))}  f^\delta (\nu ,b) db,\quad\nu \in i{\Bbb
R}_\varepsilon,\, x \in X ,
$$and thus $P^\delta$ takes $\Cal S (i\Bbb R _\varepsilon \times B) _{\check \delta ,e}$ into $\Cal S (i\Bbb R _\varepsilon \times B)
^\delta_e$.
Taking trace we see that $Q$ maps $\Cal S (i\Bbb R _\varepsilon \times B)
^\delta_e$ into $\Cal S (i\Bbb R _\varepsilon \times B) _{\check \delta ,e}$.
$\hfill \square$
\enddemo
From Proposition 2.2.6i), we see that $\Cal H (\Bbb C \times B ) ^\delta$ is a dense subspace of
$\Cal S (i\Bbb R _\varepsilon \times B)^\delta$.
Furthermore, given $F \in {\Cal S (i{\Bbb R}_\varepsilon \times B)}^\delta$, we can find
$f \in {\Cal S (i{\Bbb R}_\varepsilon \times B)}_{\check \delta,e}$ such that $F = f^\delta$.
We will use this identification in the following.
Consider $\delta \in \hat K _M$. Define the evaluation map:
$$
 (ev f) (\nu ) =  f^\delta (\nu) \equiv f^\delta (\nu, eM).
$$
We see that $\delta (m)f^\delta (\nu) = f^\delta (\nu)$,
hence $f^\delta (\nu) \in  \Hom (V_\delta ,V_\delta ^M)$.
The evaluation map is a homeomorphism between
$\Cal S (i\Bbb R _\varepsilon \times B) ^\delta$ and $\Cal S (i\Bbb R _\varepsilon , \Hom (V_\delta ,
V_\delta ^M)) \equiv \Cal S (i\Bbb R_\varepsilon) ^\delta$, and between
$\Cal H (\Bbb C \times B ) ^\delta$ and $\Cal H (\Bbb C ,\Hom (V_\delta ,V_\delta
^M)) \equiv \Cal H (\Bbb C) ^\delta$.
Lets consider $f \in  \Cal S (i\Bbb R _\varepsilon \times B) ^\delta_e$. Then:
$$\split
&\int _B e^{(-\nu + \rho ) (A(x,b))} f^\delta (- \nu ,b) db = \int _B e^{(\nu + \rho ) (A(x,b))}  f^\delta (\nu ,b) db
\Leftrightarrow\\
&\int _K e^{(-\nu + \rho ) (A(x,kM))} \delta (k) dkf^\delta (- \nu) = \int _K e^{(\nu + \rho ) (A(x,kM))} \delta (k)
dkf^\delta (\nu)
\Leftrightarrow\\
&\Phi _{-\nu , \delta } (x) f^\delta (- \nu) =\Phi _{\nu , \delta } (x) f^\delta ( \nu) ,\, \\
\endsplit\
$$where $\Phi _{\nu , \delta } (x) = \int _K e^{(\nu + \rho ) (A(x,kM))} \delta (k)dk$
is a generalized spherical function (or Eisenstein integral).
We will consider the spaces
$\Cal H (\Bbb C) ^\delta _e$ in $\Cal H (\Bbb C) ^\delta$ and $\Cal S (i\Bbb R_\varepsilon) ^\delta _e$
in $\Cal S (i\Bbb R_\varepsilon) ^\delta$ of functions
satisfying:
$$
\Phi _{-\nu , \delta } (x) f^\delta (- \nu) =\Phi _{\nu , \delta } (x) f^\delta ( \nu), \, \nu \in i{\Bbb
R}_\varepsilon,\, x \in X, \tag 1
$$
and from above we see that we have homeomorphism between
$\Cal H (\Bbb C \times B ) ^\delta _e$ and $\Cal H (\Bbb C) ^\delta _e$ and
between
$\Cal S (i\Bbb R _\varepsilon \times B) ^\delta _e$ and $\Cal S (i\Bbb R _\varepsilon ) ^\delta _e$.
Since $G/K$ is a symmetric space of rank one, we see from [He2,\,Chapter\,II,\,Corollary\,6.8] and [He1, Chapter\,V,\,Theorem\,3.5] that $\dim V_\delta ^M = 1$.
Let $v$ span $V_\delta ^M$ and let $v_1,v_2,\dots ,v_{d(\delta )}$ with $v_1 = v$ be an orthonormal basis of $V_\delta$.
\proclaim {Lemma 2.2.7} Let
$$
\varphi _{\nu , \delta } (x)= \langle \Phi _{\nu , \delta } (x) v,v \rangle = \int _K e^{(\nu + \rho ) (A(x,kM))} \langle \delta (k)
v,v\rangle dk,
$$and let
$$
\varphi _{\nu , \delta } ^j (x)= \langle \Phi _{\nu , \delta } (x) v,v_j\rangle = \int _K e^{(\nu + \rho ) (A(x,kM))} \langle \delta (k)
v,v_j \rangle dk,\, 1\le j \le d(\delta).
$$Then $\varphi _{\nu , \delta } ^j (ka \cdot x_o) = \langle \delta (k) v,v_j \rangle \varphi _{\nu , \delta } (a \cdot
x_o)$, where $x_o = eK,\,k\in K,\,a\in A$.
\endproclaim
\demo {Proof}
Define $F: X \to V_\delta $ by:
$$
F(x) =  \int _K e^{(\nu + \rho ) (A(x,kM))} \delta (k) v dk = \Phi _{\nu , \delta } (x) v .
$$Then $\varphi _{\nu , \delta } ^j (x) = \langle F(x) , v_j \rangle$ and since $F(k \cdot x) = \delta (k) F(x)$ we have
$F(a\cdot x_o) \in V_\delta ^M$ ($\delta (m) F(a \cdot x_o) = F(am \cdot x_o) = F(a\cdot x_o),\,m\in M = Z_K (A)$). Since
$\dim V_\delta ^M = 1$ we deduce: $F(a\cdot x_o) = \varphi _{\nu , \delta } (a \cdot x_o)v$. Then:
$$
\varphi _{\nu , \delta } ^j (ka \cdot x_o) = \langle F( ka\cdot x_o) , v_j \rangle = \langle \delta (k) F( a\cdot x_o) , v_j \rangle
= \langle \delta (k) v,v_j \rangle \varphi _{\nu , \delta } (a \cdot x_o) .
$$
$\hfill \square$
\enddemo
Since $f^\delta (\nu) \in \Hom (V_\delta ,V_\delta ^M))$ we see from Lemma 2.2.7 that $(1)$ is equivalent to:
$$
\varphi _{-\nu , \delta } (a\cdot x_o) f^\delta (-\nu) =
\varphi _{\nu , \delta } (a\cdot x_o) f^\delta (\nu),\quad \nu \in i{\Bbb R}_\varepsilon .\tag 2
$$
We can determine $\varphi _{\nu , \delta }$ quite explicitly in terms of the hypergeometric
functions (See [He2,\,Chapter\,III,\,Theorem 11.2]).
As a corollary we get:
\proclaim {Lemma 2.2.8} The functions $\varphi _{\nu , \delta }$ satisfies the symmetry condition:
$$
\varphi _{\nu , \delta } (a\cdot x_o) = \varphi _{-\nu , \delta } (a\cdot x_o) \frac {p_\delta (\nu)}{p_\delta (-\nu )},
$$
where ${p_\delta (\nu)}$ is a polynomial of the form:
$$
{p_\delta (\nu)} = ( \nu +\rho +s -1) \cdots (\nu +\rho),\, s\in \Bbb N \quad
\text {or} \quad {p_\delta (\nu)} \equiv 1 .
$$
\endproclaim
\demo {Proof}
[He2,\,Chapter\,III,\,Corollary\,11.3].
$\hfill \square$
\enddemo
So $(2)$ is equivalent to:
$$
{p_\delta (-\nu)}f^\delta (-\nu)={p_\delta (\nu)}f^\delta (\nu),\, \nu \in i{\Bbb
R}_\varepsilon .  \tag 3
$$Consider the set of even functions in $\Cal H (\Bbb C) ^\delta$, $\Cal H (\Bbb C) ^\delta _1$ and the set of even functions in
$\Cal S (i{\Bbb R}_\varepsilon) ^\delta$, $\Cal S (i{\Bbb R}_\varepsilon)^\delta _1$. Then we have the following lemma:
\proclaim {Lemma 2.2.9}
The map:
$$
G(\nu) \to  F(\nu )={p_\delta (-\nu)}G(\nu)
$$is a homeomorphism of $\Cal S (i\Bbb R
_\varepsilon ) ^\delta _1$ onto $\Cal S (i\Bbb R
_\varepsilon ) ^\delta _e$ taking $\Cal H (\Bbb C) ^\delta _1$ to
$\Cal H (\Bbb C ) ^\delta _e$.
\endproclaim
\demo {Proof}
The map clearly is continuous, and it takes
$\Cal H (\Bbb C ) ^\delta _1$ into $\Cal H (\Bbb C ) ^\delta _e$ and ${\Cal S (i{\Bbb R}_\varepsilon)}
^\delta _1$ into ${\Cal S (i{\Bbb R}_\varepsilon)} ^\delta _e$.
Since $p_\delta$ is a nonzero polynomial we see that the map is injective.
Let $F \in \Cal S (i\Bbb R _\varepsilon ) ^\delta _e$, and consider the function:
$$
G(\nu) = \frac {F(\nu)}{p_\delta (-\nu)}.
$$
$G$ is even:
$$
G(-\nu) = \frac {F(-\nu)}{p_m (\nu)} \cdot \frac {p_m (-\nu)}{p_m (-\nu)}
= \frac {F(\nu)}{p_m (\nu)} \cdot \frac {p_m (\nu)}{p_m (-\nu)}
=\frac {F(\nu)}{p_m (-\nu)} =G(\nu).
$$
Since $\nu$ and $-\nu$ are not both roots for $p_\delta$ we see that one of the expressions
$G(\nu) =\frac {F(\nu)}{p_m (-\nu)} = \frac {F(-\nu)}{p_m (\nu)}$ always will be welldefined,
and then $G$ will satisfy the same differentiability and growth conditions as $F$.
Hence the map is surjective, and by the closed graph
Theorem the map is a homeomorphism with the required properties.
$\hfill \square$
\enddemo
Via the classical Fourier transform it is easy to verify that
$\Cal H (\Bbb C) ^\delta _1$ is dense in $\Cal S (i{\Bbb R}_\varepsilon) ^\delta _1$,
which yields that $\Cal H (\Bbb C ) ^\delta _e$ is dense in
$\Cal S (i\Bbb R _\varepsilon ) ^\delta _e$, and thus we get the desired result:
\proclaim {Theorem 2.2.10} Let $\v \ge 0$ and let $\delta \in \hat K _M$, then:
\roster
\item "(i)" $\Cal H (\Bbb C ) ^\delta _e$ is dense in $\Cal S (i\Bbb R
_\varepsilon ) ^\delta _e$.
\item "(ii)" $\Cal H (\Bbb C \times B) _{\delta  ,e}$ is a dense subspace of
${\Cal S (i{\Bbb R}_\varepsilon \times B)}_{\delta ,e}$.
\item "(iii)" $\Cal H _e (\Bbb C \times B)$ is a dense subspace of ${\Cal S (i{\Bbb R}_\varepsilon \times B)}_e$.
\endroster\endproclaim
\demo {Proof}
(ii) Since $\Cal H (\Bbb C \times B) _{\check \delta ,e}$
is dense in $\Cal S (i\Bbb R _\varepsilon \times B) _{\check \delta ,e}$
for all $ \delta \in \hat K _M$.
\flushpar
(iii) We have
$f \in {\Cal S (i{\Bbb R}_\varepsilon \times B)}_{e}
\Leftrightarrow f_\delta \in {\Cal S (i{\Bbb R}_\varepsilon \times B)}_{\delta ,e}
\, \forall \delta \in \hat K _M$ .
\flushpar
$\Leftarrow:\,{\Cal S (i{\Bbb R}_\varepsilon \times B)}_{e}$ is closed in ${\Cal S (i{\Bbb R}_\varepsilon \times
B)}$.
\flushpar
$\Rightarrow:\,$In (SC2) replace $x$ by $k^{-1} \cdot x$, use $A(k^{-1} \cdot x ,b) = A(x, k\cdot b)$, make the substitution
$b \mapsto k^{-1} \cdot
b$, multiply by $\chi_\delta (k^{-1})$  and integrate over $K$. The result is:
$$
\int _B e^{(-\nu + \rho ) (A(x,b))} f_{\delta} (- \nu ,b) db = \int _B e^{(\nu + \rho ) (A(x,b))}  f_{\delta} (\nu ,b) db.
$$
which exactly means that
$f_\delta \in {\Cal S (i{\Bbb R}_\varepsilon \times B)}_{\delta ,e}$.

By Theorem 2.2.10ii), every finite sum of the type $\sum _{\delta \in \hat K _M} f_\delta,\, f_\delta \in
{\Cal S (i{\Bbb R}_\varepsilon \times B)}_{\delta,e}$, can be approximated by a finite sum $\sum _{\delta \in \hat K _M} g_\delta,\,
g_\delta \in \Cal H  (\Bbb C \times B) _{\delta,e}$, and the theorem then follows by Theorem 2.2.5.
$\hfill \square$
\enddemo
{\bf Remark:} As the primary goal was to show Theorem 1 via Anker's method, Theorem\,1.2.10iii) is very much an important result.
But since we can only handle the $K$-finite case, the important result is actually Theorem\,1.2.10i).

Later we will need some estimates of the matrix coefficients of the generalized spherical functions (as for spherical functions
in [An]).
Let $P=MAN$ be the minimal parabolic subgroup of $G$. Consider the characters $\{ma_tn \mapsto e^{\nu t},\, \nu \in \Bbb C \}$ of $P$.
The spherical principal series representations ${\pi}_\nu$ on $G$ are then the induced representations coming from these characters. They are
realized on $L^2(K/M)$ via the formula:
$$
\{ {\pi}_\nu (g) f \} (k) = e^{-(\nu + \rho) H(g^{-1} k)} f (\kappa (g^{-1} k)),\quad f\in L^2(K/M).
$$
The restriction of ${\pi}_\nu$ to $K$ is: $\{ {\pi}_\nu (k_1) f \} (k_2)=f ({k_1}^{-1} k_2)$.
\proclaim {Proposition 2.2.11} Let $\p \in C^\infty (K/M)$, and consider the function:
$$ \psi _\nu (g) =  \int _K \p(k) e^{- (\nu +\rho ) (H(g^{-1}k))} dk.
$$Given $D,E \in U(\goth g)$, there is a constant $c$ and elements $D _l \in U(\goth
k),\,1\le l \le M$
such that:
$$
| \psi _\nu (D;g;E) | < c(1 + |\nu |) ^{deg D +deg E}\varphi _{Re \nu} (g) \{\sum _{l=1} ^M \sup _{k \in
K} | \p(k;D _l)| \} .
$$
\endproclaim
\demo {Proof}
We can write $| \psi _\nu (D;g;E) |$ as:
$$\split
| \int _K  \p(k) \pi _\nu (D;g;E)  1_{K/M} dk |
&= | \int _K  \p(k) \pi _\nu (D) \pi _\nu (g) \pi _\nu (E) 1_{K/M} dk | \\
&= |\langle \pi _\nu (g) \pi _\nu (E) 1_{K/M} , \pi _{\bar \nu} (D^* ) \overline {\p(k)} \rangle _2 | \\
&\le \| \pi _\nu (g) \pi _\nu (E) 1_{K/M} \| _1 \| \pi _{\bar \nu} (D^* )  \overline {\p} \| _\infty ,
\endsplit\
$$
by the H\"older inequality. $\pi _\nu (E) 1_{K/M}$ can be considered as a function on $K/M:\,\xi$, hence we have:
$$
[\pi _\nu (g) \xi] (k) =  e^{- (\nu +\rho ) H(g^{-1}k)} \xi (\kappa(g^{-1}k)).
$$
Again using H\"older this gives us:
$$
\| \pi _\nu (g) \pi _\nu (E) 1_{K/M} \| _1 \le \|\xi (k) \|_\infty \| \pi _\nu (g) 1_{K/M} \| _1 =
\| \pi _\nu (E) 1_{K/M} \| _\infty \varphi _{Re \nu} (g).
$$We are left by estimating the two norms:
\roster
\item "(i)" $\| \pi _\nu (E) 1_{K/M} \| _\infty$
\item "(ii)" $\| \pi _{\bar \nu} (D^* )  \p(k) \| _\infty$
\endroster

From [An, p.336] we get an estimate on (i) $\| \pi _\nu (E) 1_{K/M} \| _\infty \le c (1 + |\nu |) ^{deg E}$.
So consider (ii). Introduce the auxiliary function $f_\nu (g) = e^{- (\nu +\rho ) H(g)}$, then by the Leibniz rule of
differentiation:
$$\split
\pi _{ \nu} (D^* )   \p(k) &= \sum_{deg D =deg D^\prime +deg D^{\prime \prime}}f_\nu (D^\prime ;k)  \p (\kappa (D^{\prime \prime};
k))\\ &=
\sum_{deg D =deg D^\prime +deg D^{\prime \prime}} f_\nu (D^\prime;k)  \p (\kappa (k;Ad (k^{-1})D^{\prime \prime})).
\endsplit
$$
We see that $f_\nu ( D^\prime;k) = \{\pi _{ \nu} ({D^\prime}^* ) 1_{K/M}\} (k)$.
There exists functions $\xi _1,\dots ,\xi _m$ in $C^\infty (K)$ and elements $D_1,\dots , D_m$ in $U(\goth g)$ of
degree $\le \deg D$ such that $Ad (k^{-1}D^{\prime \prime}) = \sum _{l=1} ^m \xi _l (k) D_l$. By the PBW-Theorem we can write:
$$
U(\goth g) = U(\goth k) \oplus U(\goth g) (\goth a +\goth n ).
$$
Let $D^{\prime } _l$ be the projection of $D_l$ on $ U(\goth k)$. Then:
$$\split
| \p (\kappa (k;Ad (k^{-1})D^{\prime \prime}))| &= | \p (\kappa (k;\sum _{l=1} ^m \xi _l (k) D_l))|
=|\sum _{l=1} ^m \xi _l (k) \p (\kappa (k;D_l))|\\
&=|\sum _{l=1} ^m \xi _l (k) \p (\kappa (k;D^{\prime} _l))| \le c \sum _{l=1} ^m \sup _{k \in K} | \p (k;D^{\prime} _l)| ,\endsplit\
$$
and thus:
$$
\| \pi _{\bar \nu} (D^* )   \p \| _\infty \le c (1 + |\nu |) ^{deg D} \sum _{deg D^{\prime \prime} \le deg D} \sum _{l=1} ^m
\sup _{k \in K} | \p (k;D^{\prime}  _l)|.
$$
$\hfill \square$
\enddemo

\flushpar
{\nilsfont 2.3. The isomorphism of the Fourier transform on $K$-finite elements in the $L^p$-Schwartz spaces.}

In this section we show the generalization of [An] in the rank $1$ case. Let $f \in \Cal S ^p (G/K),$
$0<p \le2$.
Since $C ^\infty _c (G/K) \subset \Cal S ^p (G/K) \subset L^2 (G/K)$,
we can define the Fourier transform $\hat f \in L^2 (i\Bbb R _+ \times B,c{|c(i\nu)|}^{-2})$.
From the definition of $\Cal S ^p (G/K)$,
and the Cartan decomposition of the measure on $X$, we see that the
extension of the Fourier transform from $C ^\infty _c (G/K)$ to $\Cal S ^p (G/K)$ is trivial, that is:
$$
\Cal H f(\nu,b) = \hat f (\nu,b) = \int _X f(x) e^{(-\nu + \rho ) (A(x,b))} dx,
$$
for all $\nu \in i\R , \, b \in B$.
Furthermore we see that the above integral is welldefined for $(\nu ,b) \in i{\Bbb R}_\varepsilon\times B$.
We actually have:
\proclaim {Lemma 2.3.1}
Let $f \in \Cal S ^p ( G /K)$, and let $\varepsilon = \frac 2p -1$. Then
$\hat f \in C^\infty (i{\Bbb R}_\varepsilon\times B)$,
and $\hat f (\cdot ,b)$ is holomorphic in $i{\Bbb R}^o_\varepsilon$ for fixed $b \in B$.
\endproclaim
\demo {Proof}
Let
$$
\hat f (\nu,b) = \int _X f(x) e^{(-\nu + \rho ) (A(x,b))} dx.
$$
By the Cartan decomposition we have:
$$
\hat f (\nu,b) = \int _0 ^\infty \int _K f(ka_t) e^{(-\nu + \rho ) (A(k a_t x_o,b))} dk \sinh ^ n t  dt.
$$
Consider the above integral over $K$:
$$ \split
\left | \int _K f(ka_t) e^{(-\nu + \rho ) (A(k a_t x_o,b))} dk \right | &\le
\{ \sup _{k \in K} | f(k a_t)|\} \int _K  e^{(- Re \nu + \rho ) (A(k a_t x_o,b))} dk \\
&= \{ \sup _{k \in K} | f(k a_t)|\} \varphi _{Re \nu} (a_t) .\endsplit\
$$
We then see:
\flushpar
1) $e^{-2(\frac 1p -1) \rho t} {(1+t)}^{-1}\varphi _{Re \nu} (a_t) $ is a bounded function for $\nu \in i{\Bbb R}^o_\varepsilon$
(by Lemma 2.1.9).
\flushpar
2) $t \mapsto \{ \sup _{k \in K} | f(k a_t)| \} \in C({\overline \R} _+)$.
\flushpar
3) $\sup _{t>0} {( 1 + t )}^N {\varphi _o (a_t)}^{-\frac 2p} \{ \sup _{k \in K} | f(k a_t) |\} \le \sup _{g \in G} {( 1 + |g| )}^N {\varphi _o (g)}^{-\frac 2p} | f(g) | < \infty $.
\flushpar
4) $ t \mapsto e^{2(\frac 1p - 1 )\rho t } {(1+t)}
\{ \sup _{k \in K} | f(k a_t)|\} \in L^1 (\Bbb R _+ , \sinh ^n t  dt )$
(Lemma 2.1.9).

We thus see that:
\flushpar
5) $\hat f (\nu , b)$ is welldefined for all $\nu \in i{\Bbb R}^o_\varepsilon ,\, b \in B$.
\flushpar
6) By Morera's Theorem we see that $\hat f (\cdot , b) $ is holomorphic in $i{\Bbb R}^o_\varepsilon$.
Differentiability follows from Lebesque's Dominated Convergence Theorem.
$\hfill\square$
\enddemo

\proclaim {Theorem 2.3.2}
Let $0<p \le 2, \, \varepsilon = \frac 2p -1$.
The Fourier transform $\Cal H$ is an injective and continuous homomorphism from $\Cal S ^p( G /K)$ into ${\Cal S (i{\Bbb R}_\varepsilon \times B)}_e$.
\endproclaim
\demo {Proof} Let $f \in \Cal S ^p(G /K)$.
By Lemma 2.3.1 we see that $\Cal H f (\cdot , b) $ is holomorphic in $i{\Bbb R}^o_\varepsilon$ for fixed $b \in  B$.
Now we observe that ${\Cal S (i{\Bbb R}_\varepsilon \times B)}_e$ and it's topology also is determined by the set of seminorms:
$$
\tilde \tau _{P,M,N} ^\varepsilon (f) = \sup _{\nu \in i{\Bbb R}_\varepsilon ,\, k \in K}
\left | \left ( P\left (\frac {\partial}{\partial \nu} \right) ,{\Omega _K}^M  \right ) \{ {({\nu}^2 -\rho ^2 +d)}^N f(\nu,k)\}\right | ,
$$
where $P$, $M$ and $N$ are as before, and $d$ is a constant such that ${\nu}^2 -\rho ^2 +d \ne 0$.
Fix $\tilde \tau _{P,M,N} ^\varepsilon $. Then
$$\split
&\left ( P\left (\frac {\partial}{\partial \nu} \right ) ,{\Omega _K}^M  \right ) {({\nu}^2 -\rho ^2 +d)}^N \Cal H f (\nu
,k)\\ &=
\left ( P\left (\frac {\partial}{\partial \nu} \right ) ,{\Omega _K}^M  \right ) {({\nu}^2 -\rho ^2 +d)}^N  \int _X f(x) e^{(-\nu + \rho )
(A(x,kM))} dx \\
&= \int _X ( {\Omega _K}^M {(\Delta +d)}^N f)(x) P\left (\frac {\partial}{\partial \nu} \right )e^{(-\nu + \rho )
(A(x,kM))} dx , \endsplit\
$$
where $ {\Omega _K}^M f (x) = {\Omega _K}^M f (g) = f({\Omega _K}^M ;g),\,x=gK$. Now:
$$\split
\left | P \left (\frac {\partial}{\partial \nu} \right ) e^{(-\nu + \rho ) (A(x,kM))} \right |  &=
\left | P (-A(x,kM)) e^{(-\nu + \rho ) (A(x,kM))} \right | \\
&\le c {( 1 + |x| )}^{deg P} e^{(- Re \nu + \rho ) (A(x,kM))} .\endsplit\
$$
As before (Lemma 2.3.1), we use the Cartan decomposition and the estimate $\sinh ^n t \le e^{2\rho t}$:
$$\split
&\int _X \left | ( {\Omega _K}^M {(\Delta +d)}^N f)(x) P\left (\frac {\partial}{\partial \nu} \right )e^{(-\nu + \rho )
(A(x,kM))} \right |  dx \le \\
& c
\int _0 ^\infty \{ \sup _{k^\prime \in K} |( {\Omega _K}^M {(\Delta +d)}^N f) (k^\prime a_t)| \} {( 1 + t )}^{deg P} \int _{K} e^{(- Re \nu + \rho )
(A(k^\prime a_t,kM))} dk^\prime e^{2\rho t} dt \\
& =c\int _0 ^\infty \{ \sup _{k \in K} |( {\Omega _K}^M {(\Delta +d)}^N f) (k a_t)|\} {( 1 + t )}^{deg P} \varphi _{Re \nu } (a_t) e^{2\rho t} dt
.\endsplit\
$$
By Lemma 2.1.9 ii) there exists a constant $c$ such that
$|\varphi _{Re\nu} (a_t) | \le c(1 +t )e^{(\frac 2p - 2) \rho
t},$ for $t>0$.
Let
$$
\sigma (f)  = \sup _{g \in G} {(1+ |g|)}^{M'} {\varphi _o (g)}^{- \frac 2p} |({\Omega _K}^M {(\Delta +d)}^N f)(g)|,
$$
where we choose $ M' \in \Bbb N$ s.t. $ M' > deg P +1 +\frac 2p $. Then from Lemma 2.1.9 i) we get:
$$
|({\Omega _K}^M {(\Delta +d)}^N f) (g)| \le c \sigma (f) {(1+| g |)}^{- M' -1 +\frac 2p} e^{-\frac 2p \rho | g | }.
$$
Putting this together we get:
$$
\tilde \tau _{P,M,N} ^\varepsilon (\Cal H f ) \le c \sigma (f)
$$and $\sigma$ is obviously bounded by a sum of seminorms for $\Cal S ^p( G /K)$,
and we have shown continuity.

Injectivity follows from Theorem 2.1.3.
$\hfill\square$
\enddemo
Note, that the restriction to powers of $\Omega _K$ is not crucial,
we could prove the theorem for the
topology with $\Omega _K$ replaced by $E \in U (\goth g)$ in the same fashion.

Now, of course, we would like to show surjectivity of $\Cal H$. In order to generalize the proof from [An], we
need a cut-off function, $\omega _j$, w.r.t. to the first variable of $H (\cdot, \cdot)$, in order to control
the involved functions in the commutative diagram below:
$$
\diagram
& {h \in {\Cal H}_e (\Bbb C \times B) }   \\
 {f  \in C_c^\infty ( G /K)} \urto^{\Cal H}_\cong \rrto_{\Cal R}^\cong& &
{H \in \Cal R C_c^\infty ( G /K)}  \ulto_{\Cal F}^\cong \\
\enddiagram
$$Unfortunately the symmetry conditions imposed on $h$ and $H$ are rather obscure,
as we have already encountered,
therefore looking at the function $H$ decomposed w.r.t. to the cut-off function,
$$
H (t,b) = \omega _j (t) H (t,b) +(1-\omega _j (t) ) H (t,b),\, H_j(t,b) = (1-\omega _j (t) ) H (t,b),
$$
will not necessarily give us a function $H_j$ in $\Cal R C_c^\infty ( G /K)$.
This means that we will be unable to define functions
$h_j = \Cal F H_j \in {\Cal H}_e (\Bbb C \times B)$ and $f_j = \Cal R ^{-1} H_j \in C_c^\infty ( G /K)$,
as e.g.,\,$f_j$, though welldefined, doesn't need to have compact support.
A way to avoid the difficulties created by the symmetry conditions is to look at fixed $K$-types,
as we did with the density argument.

So consider the left-regular representation of $K$ ($l$) on the Fr\'echet space $\Cal S ^p( G /K)$.
It is straightforward to check that this representation is a smooth Fr\'echet representation of $K$.
With $\delta \in \hat K _M$ acting on $V_\delta$,
consider the spaces $ C_c^\infty ( G /K ,\Hom (V_\delta ,V_\delta ))$
and $\Cal S ^p (G /K,\Hom (V_\delta ,V_\delta ))$ of compactly supported differentiable
functions on $G/K$, respectively $L^p$-Schwartz functions on $G /K$,
taking values in $\Hom (V_\delta ,V_\delta)$ (looking at the marix-entries).
Define $\Cal S ^p (G /K) ^\delta \equiv \{F \in \Cal S ^p (G /K ,\Hom (V_\delta ,V_\delta ))
: F( k\cdot x ) = \delta (k) F (x) \}$,
and define
$C_c^\infty ( G /K) ^\delta =\Cal S ^p( G /K) ^\delta \cap
C_c^\infty ( G /K ,\Hom (V_\delta ,V_\delta))$.
For $f \in \Cal S ^p( G /K)$ we define:
$$
'P ^\delta f (x ) = f^\delta (x ) = d (\delta) \int _K \delta (k) f ( k^{-1} \cdot x) dk.
$$
We see that $'P^\delta$ takes $ C_c^\infty ( G /K)$ into $C_c^\infty ( G /K) ^\delta$, and
$\Cal S ^p (G /K)$ into $\Cal S ^p (G /K) ^\delta$.
Let furthermore $C_c^\infty ( G /K) _{ {\delta}}$ and $\Cal S ^p (G /K) _{ {\delta}} $ denote the space of
$K$-finite functions in $C_c^\infty ( G /K)$, respectively in $\Cal S ^p (G /K)$, of type $\delta$.
The continuous projection of $C_c^\infty ( G /K)$, respectively of $\Cal S ^p (G /K)$,
onto $C_c^\infty ( G /K) _{ {\delta}}$ and $\Cal S ^p (G /K) _{ {\delta}}$
is given by:
$$
'P_{\delta} f (x ) = f_\delta (x ) = d (\delta) \int _K \chi _\delta (k^{-1}) f ( k^{-1} \cdot x) dk.
$$
\proclaim {Proposition 2.3.3}
Let $\delta \in \hat K _M$.
\roster
\item "(i)" The map:
$$
Q:\, F(x) \to Tr (F(x))
$$is a homeomorphism of $\Cal S ^p (G/K) ^\delta$ onto
${\Cal S ^p(G/K)}_{\check \delta}$, with inverse
$f \to {'P} ^\delta f = f^\delta$. Furthermore $Q$ takes $C_c^\infty ( G /K) ^\delta$ onto $C_c^\infty ( G /K) _{\check
\delta}$. Here $Tr$ denotes trace in $\Hom (V_\delta ,V_\delta )$.
\item "(ii)" The maps:
$$\gather
'P_{ \delta}:\,{\Cal S ^p(G/K)} \to {\Cal S ^p(G/K)}_{ \delta}, \quad \text {and}\\
'P ^\delta:\,{\Cal S ^p(G/K)} \to {\Cal S ^p(G/K)}^{ \delta},\endgather
$$are continuous open surjections, and the images are closed in ${\Cal S ^p(G/K)}$, respectively in
$\Cal S ^p (G /K,\Hom (V_\delta ,V_\delta ))$.
\endroster\endproclaim
\demo {Proof}
As for proposition 2.2.6.
$\hfill \square$
\enddemo
\proclaim {Definition 2.3.4}
Let $\delta \in \hat K _M$. For $f \in \Cal S ^p (G/K) _{\check \delta}$
the $\delta$-spherical transform is defined
by:
$$
\split
\Cal H ^\delta f (\nu) &=  (ev \circ P ^\delta \Cal H f) (\nu )\\
&= d(\delta) \int _K \int _X f( x)  e^{(-\nu + \rho ) (A(x,k M))} dx \delta (k^{-1}) dk \\
&= d(\delta) \int _X f(x) \int _K e^{(-\nu + \rho ) (A(x,kM))} \delta (k^{-1}) dk dx\\
&= d(\delta) \int _X f(x) \Phi _{- \bar \nu , \delta} (x) ^* dx ,\\
\endsplit\
$$
where: $\Phi _{- \bar \nu , \delta} (x) ^* = \int _K e^{(-\nu + \rho ) (A(x,kM))} \delta (k^{-1}) dk$
is the adjoint of the generalized spherical function $\Phi _{- \bar {\nu} , \delta}$,
${}^*$ denoting the adjoint in $\Hom (V_\delta ,V_\delta )$,
and $ev$ is the evaluation map, $ev f (\nu) = f(\nu, eM)$.
\endproclaim
For the trivial representation we have
$\Phi _{- \bar {\nu} , 1 } (x) ^* =\varphi _\nu (x)$, that is, the 1-spherical transform
$\Cal H ^1$ is the "classical" spherical transform $\Cal H$.
For the $\delta$-spherical transform we have the following Paley-Wiener Theorem:
\proclaim {Theorem 2.3.5}
The $\delta$-spherical transform $f \mapsto  \Cal H ^\delta f$ is a bijection of
$C_c^\infty ( G /K) _{\check \delta}$ onto $\Cal H (\Bbb C) ^\delta _e$.
\endproclaim
\demo {Proof}
It is evident from the definition that $\Cal H ^\delta$ maps
$C_c^\infty ( G /K) _{\check \delta}$ into $\Cal H (\Bbb C) ^\delta _e$,
and also that $\Cal H ^\delta$ is injective since $\Cal H$ is an injective map.
For surjectivity, let $\psi \in \Cal H (\Bbb C) ^\delta _e$. The function
$\Psi (\nu , kM) \equiv Tr (\delta (k) \psi (\nu))$ is clearly a holomorphic function
of uniform exponential type on $\Bbb C \times B$. By Proposition 2.2.6 and the
discussion afterwards, we see that $\Psi \in \Cal H (\Bbb C \times B) _{\check \delta}$.
By the Paley-Wiener Theorem 2.1.5,
there exists a unique $F \in C^\infty _c (G/K)$ such that $\Psi = \Cal H F$.
By Proposition 2.3.3 the function:
$'P_{ \check \delta} F = F _{\check \delta}$ belongs to $C^\infty _c (G/K)_ {\check \delta}$,
and we have:
$$\split
\Cal H ^\delta F _{\check \delta} (\nu) &= \{ P ^\delta \Cal H {}'P_{{\check \delta}} F\} (\nu , eM) =
\{ P ^\delta P_{ {\check \delta}}  \Cal H F\} (\nu , eM) = \{ P ^\delta P_{{\check \delta}} \Psi \} (\nu , eM) \\
&=d^2 (\delta ) \int _K \int _K \Psi (\nu , u^{-1} k M) \chi _\delta (u) \delta (k^{-1}) dkdu \\
&=d^2 (\delta ) \int _K \int _K Tr (\delta (k) \psi (\nu))  \chi _\delta (u) \delta (k^{-1} u^{-1}) dkdu\\
&=d^2 (\delta ) \int _K \int _K Tr (\delta (k) \psi (\nu))  \chi _\delta (u) \delta (k^{-1})\delta (u^{-1}) dkdu\\
&=d^2 (\delta ) \int _K \delta (k^{-1})Tr (\delta (k) \psi (\nu)) dk \int _K \chi _\delta (u) \delta (u^{-1}) du\\
&=d (\delta ) \int _K \delta (k^{-1})Tr (\delta (k) \psi (\nu)) dk\\
&=  P ^\delta Tr (\delta (k) \psi (\nu))\\
&= \psi (\nu),
\endsplit\
$$
since the orthogonality relations yields:
$$
 d(\delta) \int _K \bar {\delta} _{i,j} (u) \delta (u) du = E_ {i,j}.
$$
$\hfill \square$
\enddemo

We also have the following inversion Theorem:
\proclaim {Theorem 2.3.6}
The $\delta$-spherical transform is inverted by:
$$
f(x) = c Tr \left \{ \int _0 ^\infty \Phi _{i\nu , \delta} (x) (\Cal H ^\delta f) (i\nu) |c (i\nu)| ^{-2} d\nu
\right \},\quad f \in C_c^\infty ( G /K) _{\check
\delta},
$$where
$$
 \Phi _{i\nu , \delta} (x)  = \int _K e^{(i\nu + \rho ) (A(x,kM))} \delta (k) dk.
$$
\endproclaim
\demo {Proof}
If $f \in C_c^\infty ( G /K) _{\check \delta}$, then:
$$
\split
\Cal H f (\nu , kM) &= \Cal H {}'P_{{\check \delta}} f (\nu , kM) = P_{{\check \delta}} \Cal H  f (\nu , kM)\\
&=Tr P ^\delta P_{{\check \delta}} \Cal H  f (\nu , kM) = Tr P ^\delta \Cal H {}'P_{{\check \delta}} f (\nu , kM)\\
&=Tr \Cal H ^\delta f (\nu , kM) = Tr (\delta (k)\Cal H ^\delta f (\nu)).
\endsplit
$$
So by the "classical" inversion Theorem, Theorem 2.1.2, we get:
$$
\split
f(x) &= c \int _{\Bbb R _+ \times K}  e^{(i\nu + \rho ) (A(x,kM))} Tr (\delta (k) \Cal H ^\delta f (i\nu))
  |c (i\nu)| ^{-2} d\nu \\
&= c Tr \left \{\int _{\Bbb R _+ \times K}  e^{(i\nu + \rho ) (A(x,kM))}  \delta (k) \Cal H ^\delta f (i\nu)
  |c (i\nu)| ^{-2} d\nu\right \} .
\endsplit\
$$
$\hfill \square$
\enddemo
Let $v$ span $V_\delta ^M$, and let $v_1,v_2,\dots ,v_{d(\delta )}$,
with $v_1 = v$, be an orthonormal basis of $V_\delta$.
Then the inverse $\delta$-spherical transform can be written as:
$$\split
f(x) &=
c \int _0 ^\infty  Tr ({\Phi _{i\nu , \delta} (x)}  { \Cal H ^\delta f (i\nu)})
 |c (i\nu)| ^{-2} d\nu\\
&=c \int _0 ^\infty  \sum _{l=1} ^{d(\delta)}{\Phi _{i\nu , \delta} (x)} _{l,1} { \Cal H ^\delta f (i\nu)}_{1,l}
 |c (i\nu)| ^{-2} d\nu,
\endsplit
$$
where $i,j$ denote matrix entries.
We see that ${\Phi _{i\nu , \delta} (x)} _{l,1} = \varphi ^l _{i\nu , \delta} (x)$
(see Lemma 2.2.7).

Moreover we also have a Plancherel Theorem:
\proclaim {Theorem 2.3.7}
$$
\split
\int _X |f(x)| ^2 dx &=  \int _0 ^\infty Tr \{ (\Cal H ^\delta f) (i\nu)
(\Cal H ^\delta f) ^* (i\nu) \}  |c (i\nu)| ^{-2} d\nu\\
&= c \int _0 ^\infty  \| \{ (\Cal H ^\delta f) (i\nu) \} \| _{HS} ^2 |c (i\nu)| ^{-2}
d\nu .\endsplit
$$
\endproclaim
\demo {Proof}
As above.
$\hfill \square$
\enddemo
Consider the classical Fourier transform $\Cal F$ acting on vector valued functions, and define
$C_c^\infty (\Bbb R) ^\delta _e = \Cal F ^{-1} \Cal H (\Bbb C) ^\delta _e$, that is, functions in
$C_c^\infty (\Bbb R, \Hom (V_\delta ,V_\delta ))$ satisfying certain symmetry conditions, then
we have the following commutative diagram:
$$
\CD
C_c^\infty ( G /K) _{\check \delta} @> \Cal H ^\delta >> \Cal H (\Bbb C) ^\delta _e \\
@V\Cal R VV        @AA \Cal F A \\
C_c^\infty (\Bbb R \times B) _{\check \delta ,e}  @> ev \circ {}''P  ^\delta  >> C_c^\infty (\Bbb R) ^\delta _e
\endCD
$$
where:
$$
''P  ^\delta f(t,b) = d (\delta) \int _K \delta (k) f ( t,k^{-1} \cdot b) dk,
$$
and $\Cal R$ is the Radon transform.
We define a new transform $\Cal T$ by:
$$
\Cal T =ev \circ {}''P ^\delta \Cal R .
$$More exactly we have:
$$
\Cal T f (t)= e^{\rho t } \int _{K \times N} f (k a_t n )\delta (k^{-1}) dk dn ,
$$which can be thought of as a generalized Abel transform.
\proclaim {Theorem 2.3.8} The transform
$\Cal T$ is an isomorphism between $C_c^\infty ( G /K) _{\check \delta}$ and $C_c^\infty (\Bbb R) ^\delta _e$. Moreover
$\supp f \subset \bar B (0, R)$ if and only if $\supp \Cal T f \subset
\bar B (0, R) $.
\endproclaim
\demo {Proof}
i) By definition.
\flushpar
ii) $\Leftarrow )\, $As Proposition 2.1.6.$\,\Rightarrow )\,$Let $ g \in C_c^\infty (\Bbb R) ^\delta _e$ such that $\supp g \subset \bar B (0, R)$. Then $\Cal F g \in \Cal H (\Bbb C) ^\delta _e$,
with matrix entries $\{\Cal F g \} _{i,j}$ in $\Cal H ^R (\Bbb C)$.
Now consider the function $G (z , kM) = Tr (\delta (k) \{\Cal F g \} (z) )$ in
$\Cal H _e ^R (\Bbb C \times B)$. The proof of Theorem 2.3.6 yields that:
$$
\Cal T ^{-1} g  (x) = \Cal H ^{-1} G (x),
$$
and then by Proposition 2.1.7, we see that $\supp \Cal T  ^{-1} g \subset \bar B (0, R)$.
$\hfill \square$
\enddemo
Now we arrive at the essential theorem:
\proclaim {Theorem 2.3.9}
Let $0 <p\le 2, \varepsilon = \frac 2p -1$. Then the $\delta$-spherical transform $\Cal H ^\delta$ is a topological
isomorphism between $\Cal S ^p (G /K) _{\check \delta} $ and $\Cal S (i\Bbb R_\varepsilon) ^\delta _e$. The inverse transform is
given by Theorem 2.3.6.
\endproclaim
\demo {Proof}
a) We can write $\Cal H ^\delta$ as a composition of continuous operators:
$$
\Cal H ^\delta  =ev \circ  P ^\delta \Cal H .
$$
Hence $\Cal H ^\delta$ is a injective continuous homomorphism into $\Cal S (i\Bbb R_\varepsilon) ^\delta _e$.

b) We now want to show surjectivity. By density, Theorem 2.2.10i), it is enough to show that the inverse $\delta$-spherical transform is
continuous as a map from $\Cal H (\Bbb C) ^\delta _e$ to $C_c^\infty ( G /K) _{\check \delta}$,
with the topologies induced by
$\Cal S (i\Bbb R_\varepsilon) ^\delta _e$ and $\Cal S ^p (G /K) _{ {\delta}}$.
So let $f  \in C_c^\infty (G /K) _{\check \delta} $, $h = \Cal H ^\delta f$
and $H = \Cal T f$ as in the commuting diagram below.
$$
\diagram
& {h \in {\Cal H} (\Bbb C)_e ^\delta }   \\
 {f  \in C_c^\infty ( G /K) _{\check \delta}} \urto^{\Cal H ^\delta}_\cong \rrto_{\Cal T}^\cong& &
{H \in C_c^\infty (\Bbb R) ^\delta _e}  \ulto_{\Cal F}^\cong \\
\enddiagram
$$
Let $\sigma _{D,E,N} ^p $ be a seminorm for $\Cal S ^p( G /K)$:
$$
\sigma _{D,E,N} ^p  (f) = \sup _{g \in G} {(1+ |g|)}^N {\varphi _o (g)}^{- \frac 2p} |f(D;g;E)|.
$$
We will consider the function: $F(g) = {(1+ |g|)}^N {\varphi _o (g)}^{- \frac 2p} |f(D;g;E)|$.
Our goal now is to estimate $F$ on the intervals $|g| \in ]j,j+1]$.
Remark that all positive constants appearing below may depend on $N$ and $p$, but not on $f$.

Step 1 ($[0,2]$): By the inversion formula, Theorem 2.3.6, we get:
$$
f(D;g;E) = c Tr \left \{ \int _0 ^\infty  \Phi _{i\nu ,\delta}(D;g;E)
h(i\nu)  |c(i\nu)|^{-2} d \nu \right \}.
$$
Using Lemma 2.1.9 iii) and Proposition 2.2.11, on the functions $\p = \delta _{l,1}$, we get:
$$
|f(D;g;E)| \le c \varphi _o (g) \int _{{\Bbb R}_+} (1+ |\nu|)^R \sum _{l=1} ^{d(\delta)} |h_{1,l}(i\nu)| d \nu,
$$
for some $R \in \Bbb N$. We thus get:
$$
\split
\sup _{|g| \in [0,2]} F(x) &\le c \int _{{\Bbb R}_+} (1+ |\nu|)^{R}  \sum _{l=1} ^{d(\delta)} |h_{1,l}(i\nu)| d \nu\\
&\le \sup _{\nu \in \Bbb R _+} (1+ |\nu|)^{R+2}  \sum _{l=1} ^{d(\delta)} |h_{1,l}(i\nu)| \\ &= c
\sum _{l=1} ^{d(\delta)} \tau _{1,R+2} ^o (h_{1,l}) ,\endsplit\
$$
using the compactness of $[0,2]$.

Step 2:
For $j \in \Bbb N$ we introduce an even auxiliary function $\omega _j \in C_c ^\infty (\Bbb R)$ defined such that:
$$
\omega _j (t) =
\left\{
\aligned &1,\quad t\in [0,j-1[\\ &0, \quad t\in [j,\infty [ \quad ,\endaligned \right .
$$
and $\omega _{j+1} (t) = \omega _j (t-1)$. We write $H$ as:
$$
\split
H(t)= &p_\delta \left ( -\frac {\partial}{\partial t} \right )
\{ (1 - \omega _j ) (t) \Cal F ^{-1}\{ h (\cdot )
 p_\delta ^{-1} (-\cdot) \} (t)\}\\ &+p_\delta \left ( -\frac {\partial}{\partial t} \right )
 \{  \omega _j  (t) \Cal F ^{-1}\{ h (\cdot )
 p_\delta ^{-1} (-\cdot) \} (t)\} .\endsplit
$$
Consider the functions:
$$
\split
H^j (t) &= p_\delta \left ( -\frac {\partial}{\partial t} \right ) \{ (1 - \omega _j ) (t) \Cal F ^{-1}\{ h (\cdot )
 p_\delta ^{-1} (-\cdot) \} (t)\} .\\
h^j (\nu ) = {\Cal F H^j} (\nu) &= p_\delta (-\nu) \{ \Cal F \{ (1 - \omega _j ) (\cdot) \Cal F ^{-1} \{ h (\cdot )
 p_\delta ^{-1} (-\cdot) \} \} \} (\nu) .\\
\endsplit\
$$
We observe that $\{ (1 - \omega _j ) (t) \Cal F ^{-1}\{ h (\cdot )  p_\delta ^{-1} (-\cdot) \} (t)\}$
is an even function.
For $h^j$ to be in $\Cal H (\C) _e ^\delta$, it has to satisfy:
$p_\delta (-\nu) h^j (-\nu) = p_\delta (\nu) h^j (\nu)$.
$$
\split
p_\delta (-\nu) h^j (-\nu) &= p_\delta (-\nu) p_\delta (\nu) \{ \Cal F \{ (1 - \omega _j )
(\cdot) \Cal F ^{-1} \{ h (\cdot )
 p_\delta ^{-1} (-\cdot) \} \} \} (-\nu) \\
&= p_\delta (\nu) p_\delta (-\nu) \{ \Cal F \{ (1 - \omega _j ) (\cdot) \Cal F ^{-1} \{ h (\cdot )
 p_\delta ^{-1} (-\cdot) \} \} \} (\nu) \\
&= p_\delta (\nu) h^j (\nu), \endsplit\
$$
and hence $h^j \in \Cal H (\C) _e ^\delta$ and $H^j \in C_c ^\infty (\Bbb R) _e ^\delta$.
Let $f^j$ be the corresponding element of $C_c ^\infty ( G/K) _{ \check \delta}$.
Since $\omega _j$ has support in $[0,j]$,
Proposition 2.3.8 tells us that $f$ may differ from $f^j$ only inside $K \times [0,j]\times K$.

Step 3 ($]j,j+1]$): As in step 1 we get:
$$\split
|f^j (D;g;E)| &\le c \varphi _o (g) \int _{{\Bbb R}_+} (1+ |\nu|)^{R} \sum _{l=1} ^{d(\delta)} |h^j _{1,l}(i\nu)|
d \nu \\ &\le c \varphi _o (g) \sum _{l=1} ^{d(\delta)}\tau _{1,R+2} ^o (h^j _{1,l} ).\\ \endsplit\
$$
It follows that:
$$
\sup _{|g| \in ]j,j+1] } |F(g)| \le c j^N e ^{\varepsilon \rho j} \sum _{l=1} ^{d(\delta)}
\tau _{1,R+2} ^o (h^j _{1,l}),
$$
by estimates on $\varphi _\nu$ (Lemma 2.1.9).

Step 4:
We now want to find a constant $c$ and $m,n \in \Bbb N$ s.t.
$$
j^N e ^{\varepsilon \rho j} \tau _{1,R+2} ^o (h^j _{1,l}) \le
c \sum _{k=0} ^m \sup _{\nu \in i\Bbb R _\varepsilon }
{(1 +|\nu |)}^n | \nabla ^k  h_{1,l}(\nu ) |,
$$for all $l$.
In the following $\nabla$ will denote either $\nabla = \frac {d}{d\nu}$ or $\nabla = \frac {d}{dt}$.
The connection between $h^j _{1,l}$ and $H^j _{1,l}$ is:
$$
h^j _{1,l}(\nu ) = \Cal F H^j _{1,l} (\nu) =\int _{\Bbb R} H^j _{1,l}(t) e ^{-\nu t } dt,
$$
and
$$
H^j _{1,l}(t) = \Cal F ^{-1} h^j _{1,l}(t) =\frac 1{2\pi} \int _{\Bbb R} h^j _{1,l}( i\nu ) e ^{i\nu t } d\nu.
$$
Thus:
$$
\split
\tau _{1,R+2} ^o (h^j _{1,l}) &= \sup _{\nu \in i\Bbb R } (1+|\nu|)^{R+2} |h^j _{1,l}(\nu )|\\
&= \sup _{\nu \in i\Bbb R} (1+|\nu|)^{R+2} |\Cal F H^j _{1,l}(\nu )|\\
&= \sup _{\nu \in i\Bbb R } (1+|\nu|)^{R+2} |\int _{\Bbb R} H^j _{1,l}(t) e ^{-\nu t } dt| \\
&\le c \sum _{k=0} ^{R+2} \int _{\Bbb R}  \left |\nabla ^k  H^j _{1,l}(t)\right |  dt\\
&\le c \sup _{t \ge 0 } \sum _{k=0} ^{R+2}{(1+t)}^2 \left | \nabla ^k  H^j _{1,l}(t)\right|.
\endsplit\
$$
We compute the derivatives of $H^j _{1,l}(\cdot )$ by the Leibniz rule. $1 - \omega _j $
and its derivatives vanish on $[0,j-1]$, and
are bounded on $[j-1,\infty [$, uniformly in $j$. Consequently:
$$
\split
&j^N e^{\varepsilon \rho t} \tau _{1R+2} ^o (h^j _{1,l}) \\
&\le c \sum _{k=0} ^{R+2}\sup _{t\ge j -1}{(1+t)}^{N+2} e^{\varepsilon \rho t}\left |
 \nabla ^k  p_\delta \left (-\frac {\partial}{\partial t}\right )\Cal F ^{-1} \{ h _{1,l}(\cdot )  p_\delta ^{-1} (-\cdot) \} (t)\right|\\
&\le c \sum _{k=0} ^{R+2}\sup _{t\ge 0 }{(1+t)}^{N+2} e^{\varepsilon \rho t}\left | \nabla ^k  p_\delta \left(-\frac {\partial}{\partial
t} \right )\Cal F ^{-1} \{ h _{1,l}(\cdot )  p_\delta ^{-1} (-\cdot) \} (t)\right|.\endsplit\
$$
Let $P$ and $Q$ be polynomials, then by Cauchy's Theorem,
shifting integral from $i\nu$ to $i\nu +\rho \varepsilon$, we get:
$$
\split
&P(t)e^{\rho \varepsilon t }Q{\left (\frac {\partial}{\partial t} \right )\Cal F ^{-1} \{ h _{1,l}(\cdot )  p_\delta ^{-1}
(-\cdot) \} (t)}  \\
&=\frac 1{2\pi} \int _{\Bbb R} P  \left (i\frac {\partial}{\partial \nu}\right )\{Q (i\nu -\rho \varepsilon)
h _{1,l}(i\nu -\rho \varepsilon) p_\delta ^{-1} (-i\nu +\rho \varepsilon )\} e^{i\nu t } d\nu.\\
\endsplit\
$$
Hence we get:
$$ \split
& \sum _{k=0} ^{R+2}\sup _{t\ge 0 }{(1+t)}^{N+2} e^{\varepsilon \rho t}  \left | \nabla ^k p_\delta \left (-\frac {\partial}{\partial
t} \right )
\Cal F ^{-1} \{ h _{1,l}(\cdot )  p_\delta ^{-1} (-\cdot) \} (t) \right|\\
&\le  c\sum _{k=0} ^{N+2} \int _\Bbb R {(1 +|\nu|)}^{\tilde R+2} \left | \{\nabla ^k  (h_{1,l}p_\delta ^{-1}(-\cdot) )\}(i\nu -\rho \varepsilon)
 \right | d \nu ,
\endsplit\
$$
where  $\tilde R = \deg p_\delta +R$.
Now the remaining problem is to estimate $h_{1,l}(\nu) p_\delta ^{-1} (-\nu)$ and derivatives
on the boundary of the tube $i \Bbb R _\varepsilon$.
Recall that the polynomials ${p_\delta (\nu)}$ are of the form:
$$
{p_\delta (\nu)} = ( \nu +\rho +s -1) \cdots (\nu +\rho),\, s\in \Bbb N \quad \text {or} \quad {p_\delta (\nu)} \equiv 1.
$$
So we see that $p_\delta ^{-1}(-\cdot) $ has at most one singularity on the boundary.
Assume that we a singularity $\nu_o$ on the boundary of $i \Bbb R _\varepsilon$.
Since $h_{1,l} (\nu _o) = 0$, we can write $h_{1,l}$ as:
$$
\split
h_{1,l}(\nu) &= (\nu - \nu _o )  h_{1,l} ^\# (\nu) \, \text {and}\\
 h_{1,l} ^\# (\nu) &= \int _0 ^1 h_{1,l} ' (\nu _o + t (\nu - \nu _o)) dt .\\
\endsplit\
$$
We thus have: $h_{1,l} (\nu) p_\delta ^{-1}(-\nu) = p_\delta ^{-1}(-\nu)(\nu - \nu _o )  h_{1,l} ^\# (\nu)$,
where the function $p_\delta ^{-1}(-\nu)$  $(\nu - \nu _o )$ has no singularity on the boundary.
Furthermore we see that the $N$ first derivatives of $ h_{1,l} ^\#$ can be estimated by the
$(N+1)$ first derivatives of $h_{1,l}$. Hence we get:
$$
\split
& \sum _{k=0} ^{N+2} \int _{\Bbb R} {(1 +|\nu|)}^{\tilde R +2}
 | \{\nabla ^k ( h_{1,l}p_\delta ^{-1}(-\cdot) )\}(i\nu -\rho \varepsilon) | d \nu \\
&\le c\sum _{k=0} ^{N+3} \int _{\Bbb R} {(1 +|\nu|)}^{\tilde R +2}
 |\{\nabla ^k  h_{1,l}\} (i\nu -\rho \epsilon) | d \nu \\
&\le c\sum _{k=0} ^{N+3} \sup _{\nu \in i\Bbb R _\varepsilon } {(1 +|\nu|)}^{\tilde R +4}
 |\nabla ^k  h_{1,l} (\nu)|  .\\
\endsplit\
$$
All in all this gives, for some positive constant $c$
$$
\sigma _{D,E,N} ^p  (f) \le c \sum _{l=1} ^{d (\delta)} \sum _{k=0} ^{N+3} \sup _{\nu \in i\Bbb R _\varepsilon }
{(1 +|\nu|)}^{\tilde R +4}  |\nabla ^k  h_{1,l} (\nu)|.
$$
This concludes the proof.
$\hfill \square$
\enddemo
\proclaim {Corollary 2.3.10}
Let $0 <p \le 2$. The Fourier transform $\Cal H$ is a topological isomorphism between
$\Cal S ^p (G/K) _\delta$ and $\Cal S (i \Bbb R _\varepsilon \times B) _ {\delta ,e}$,
taking $C ^\infty _c (G/K) _\delta$ to $\Cal H (\Bbb C \times B)  _ {\delta ,e}$.
The inverse transform is given by Theorem 2.1.2.
\endproclaim
\demo {Proof}
Theorem 2.3.9 and the isomorphisms discussed in Proposition 2.2.6 and thereafter.
Actually let $h \in \Cal H (\Bbb C \times B) _{\check \delta ,e}$,
then we get the following estimate on $f = \Cal H ^{-1} h \in C ^\infty _c (G/K) _{\check \delta}$,
for some positive constant $c$
$$
\sigma _{D,E,N} ^p  (f) \le c {d (\delta)}^2 \sum _{l=1} ^{d (\delta)} \sum _{k=0} ^{N+3} \sup _{\nu \in i\Bbb R
_\varepsilon, k \in K }
{(1 +|\nu|)}^{\tilde R +4}  |\nabla ^k  h_{1,l} (\nu, kM)|.
$$
$\hfill \square$
\enddemo
Let $\Cal S ^p (G/K) _K$ and $\Cal S (i \Bbb R _\varepsilon
\times B) _{e,K}$ denote the $K$-finite elements of $\Cal S ^p (G/K)$ and $\Cal S (i \Bbb R _\varepsilon
\times B) _e$, then we finally get the $K$-finite version of Theorem 1:
\proclaim {Theorem 2.3.11}
Let $0 <p \le 2$. The Fourier transform $\Cal H$ is a topological isomorphism between $\Cal S ^p (G/K) _K$
and $\Cal S (i \Bbb R _\varepsilon \times B) _{e,K}$. The inverse transform is given by Theorem 2.1.2.
\endproclaim
{\bf Remark:} It is not possible via a density argument to use the above theorem to prove
Theorem 1 in general, since we don't have a polynomial or uniform bound in the estimate above for the various $K$-types.
In estimating the derivatives of $H ^j _{1,l}$, we used boundedness of finite derivations of
$\omega _j$, but unfortunately we cannot have a uniform or polynomial bound on the derivatives
of $\omega _j$, as this would imply analycity of $\omega _j$.
Hence, as there is no bound on the degree of the polynomials $p _\delta$, we
cannot have a uniform or polynomial bound on the constants $c$.

\bigskip
\flushpar
{\myfont 3. General rank.}

We will in this Section sketch how to remove the $\dim \a =1$ condition in Section 2.
So let $\goth g = \goth k \oplus \goth p$ be a Cartan
decomposition of $\goth g$, and fix a maximal subspace $\a$ in $\goth p$ ($\dim \a \ge 1$).
Denote its real dual by $\a ^*$ and its complex dual by $\a ^* _\C$.
Let $\Sigma \subset \a ^*$ be the root system of ($\goth g ,\a$),
and let $W$ be the Weyl group associated to $\Sigma$ ($\cong N_K (\a) / Z_K (\a)$).
Choose a set $\Sigma ^+ $ of positive roots, and let $\a _+ \subset \a$ and $\a ^*_+ \subset \a ^*$ be the corresponding positive Weyl chambers. We will define $\r $ as in Section 2.
Fix $\v \ge 0$, let $C^{\v \r}$ be the convex hull of the
set $W \c \v\r$ in $\a ^*$ and let $\a ^* _\v = a^* + i C^{\v \r}$
be the tube in $\a ^* _\C$ with basis $C^{\v \r}$.

In the following, the variable $\nu \in \C$ in Section 2 will correspond to the variable
$i\l \in \a _\C ^*$ and differentiation with $P \left ( \frac {\partial}{\partial \nu} \right )$
will correspond to $P \left ( \frac {\partial}{\partial \l} \right )$,
where $P \in S (\a ^*)$. On nice functions on $X = G/K$, we will define the Fourier transform as:
$$
\H f (\l,b) = \hat f (\l, b) = \int _X f(x) e^{(-i\l + \rho ) (A(x,b))} dx,
$$
for $\l \in \goth a ^* _\C,\, b \in B$, when welldefined, see [He2,\,Chapter\,III,\,\S1].
Inversion formulas, Plancherel formulas and a Paley-Wiener Theorem can be found in
[He2,\,Chapter.III].
For $\v \ge 0$, let $\Cal S (\goth a ^* _\v \times B)$ be the Schwartz space on
$\goth a ^* _\v \times B$ (replace $\R _\v ^o$ with $\goth a ^*_\v$ in definition 2.2.3),
and let $\Cal S (\goth a ^* _\v \times B) ^W$ be the subspace of $\Cal S (\goth a ^* _\v \times B)$
of functions satisfying the symmetry condition (SC):
$$
\int _B e^{(si\l + \rho ) (A(x,b))} \psi (s\l , b) db = \int _B e^{(i\l + \rho ) (A(x,b))} \psi (\l , b) db,
$$for $s \in W,\,\l \in \goth a ^* _\v$ and $ x \in X$.
For $0 <p \le 2$, $\Cal S ^p (G/K)$ denotes the $L^p$-Schwartz space on
$G/K$. The group $K$ acts naturally on $\Cal S ^p (G/K)$ and $\Cal S (\goth a ^* _\v \times B) ^W$
(on the second variable, the action preserving the symmetry condition (SC)).
Let $\Cal S ^p (G/K) _K$ and $\Cal S (\goth a ^* _\v \times B) ^W _K$
denote the $K$-finite elements of $\Cal S ^p (G/K)$
and $\Cal S (\goth a ^* _\v \times B) ^W$ respectively.
In this setup, Theorem 2.3.2 and Theorem 2.3.11 become:
\proclaim {Theorem 3.1}
Let $0 <p \le 2$ and $\v = \frac 2p - 1$.
\roster
\item "(i)" The Fourier transform
is an injective and continuous homomorphism from $\Cal S ^p (G/K)$ into $\Cal S (\goth a ^* _\v \times B) ^W$.
\item "(ii)" The Fourier transform $\Cal H$
is a topological isomorphism between $\Cal S ^p (G/K) _K$ and $\Cal S (\goth a ^* _\v \times B) ^W _K$.
\item "(iii)" The inverse transform is given by: ($\psi \in \Cal S (\goth a ^* _\v \times B) ^W _K$)
$$
\Cal H ^{-1} \psi (x) =c\int _{\goth a ^* \times B} e^{(i\l + \rho ) (A(x,b))} \hat f (\l,b) {|c(\l)|}^{-2} d\l db.
$$
\endroster\endproclaim
As in Section 1, the difficult part is to prove (ii) and (iii).
Let $C_c ^\infty (X)$ be as usual. Consider the Paley-Wiener space
$\H (\goth a ^* _\C \times B) $ (replace $\C$ with $\a ^* _\C$ and $|Re \c|$ with $|Im \c|$),
and denote by $\H (\goth a ^* _\C \times B)^W$ the subspace of $\H (\goth a ^* _\C \times B) $
of functions satisfying the symmetry condition (SC) for $\l \in \goth a ^* _\C$. Then we have:
\proclaim {Theorem 3.2} Let $0 <p \le 2$ and $\v \ge 0$. Then:
\roster
\item "(i)" $C_c ^\infty (X)$ is a dense subspace of  $\Cal S ^p (G/K)$.
\item "(ii)" $\H (\goth a ^* _\C \times B)$ is a dense subspace of $\Cal S (\goth a ^* _\v \times B) $.
\item "(iii)" $\H (\goth a ^* _\C \times B) ^W$ is a dense subspace of $\Cal S (\goth a ^* _\v \times B) ^W $.
\endroster\endproclaim
\demo {Proof}
(i) See Lemma 2.2.2.
\flushpar
(ii) See Lemma 2.2.4 and [An].
\flushpar
To prove (iii), we will consider the symmetry conditions (SC) for various $K$-types,
and again they reduce to polynomial symmetry conditions.
Recall the definitions of the various projections in Section 2,\S2,
and define the subspaces:
$\Cal S (\goth a ^* _\v \times B)_\delta ^W,\,\Cal S (\goth a ^* _\v \times B)^{\delta ,W},\,\H (\goth a ^* _\C \times B)_\delta ^W$ and $\H (\goth a ^* _\C \times B)^{\delta ,W}$
of functions of $K$-type $\delta$ satisfying (SC).
Using the evaluation map, we get isomorphisms between
$\Cal S (\goth a ^* _\v \times B)^{\delta ,W}$ and $\Cal S (\goth a ^* _\v) _W ^\delta$,
and between $\H (\goth a ^* _\C \times B)^{\delta ,W}$ and $\H (\goth a ^* _\C) _W ^\delta$, where:
$$
\split
&\H (\goth a ^* _\C) _W ^\delta \equiv
\{ F \in \H (\goth a ^* _\C, \Hom (V_\delta , V_\delta  ^M)) :
(Q ^{\check \delta} ) ^{-1} F \, \text {is $W$-invariant} \}\\
&\Cal S (\goth a ^* _\v) _W ^\delta \equiv \{ F \in \Cal S (\goth a ^* _\v,
\Hom (V_\delta , V_\delta  ^M)) : (Q ^{\check \delta} ) ^{-1} F \, \text {is $W$-invariant} \}.
\endsplit\
$$
Here $Q ^\delta (\l)$ is the Kostant $Q$-polynomial, see [He2,\,Chapter\,III,\,\S 2,3,5].
Fix an orthonormal basis $v_1,\cdots,v_{d(\delta)}$ of $V_\delta$ such
that $v_1,\cdots,v_{l(\delta)}$ span $V_\delta  ^M$.
Then the members of $\H (\goth a ^* _\C) _W ^\delta$ resp. $\Cal S (\goth a ^* _\v) _W ^\delta$ become matrix valued holomorphic functions on $\goth a ^* _\C$, respectively on $ \goth a ^* _\v$, and $Q ^{\check \delta}$ is an
$l(\delta) \times l(\delta)$ matrix whose entries are polynomials on $\goth a ^* _\C$ ([He2,\,p.287]).
Let $\H (\goth a ^* _\C) _e ^\delta$ and $\Cal S (\goth a ^* _\v) _e ^\delta$
denote the spaces of Weyl group invariants in
$\H (\goth a ^* _\C, \Hom (V_\delta , V_\delta ^M))$, respectively in
$\Cal S (\goth a ^* _\v, \Hom (V_\delta , V_\delta  ^M))$
(corresponding to $Q ^{\check \delta} \equiv I$), then we get the crucial fact:

\flushpar
{\bf Fact:}
The mapping $\psi (\l) \mapsto Q ^{\check \delta} (\l) \psi (\l)$ is a homeomorphism of
$\Cal S (\goth a ^* _\v) _e ^\delta$ onto $\Cal S (\goth a ^* _\v) _W ^\delta$,
taking $\H (\goth a ^* _\C) _e ^\delta$ onto $\H (\goth a ^* _\C) _W ^\delta$.
\flushpar
The map clearly is into, and from [He2,\,Chapter\,III,\,Lemma\,5.12]
we see that the map takes $\H (\goth a ^* _\C) _e ^\delta$ onto $\H (\goth a ^* _\C) _W ^\delta$.
We can write:
$$
Q ^{\check \delta} (\l) ^{-1} = Q _c (\l) (\det Q ^{\check \delta} (\l)) ^{-1},
$$
where $Q _c$ is a matrix whose entries are polynomials on $\goth a ^* _\C$,
and $\det Q ^{\check \delta} (\l)$ is a
product of polynomials coming from the rank one reduction, see [He2,\,p.263\,(50)] and [He2,\,Chapter\,III,\,Theorem\,4.2].
From [He2,\,\S11], we conclude that $\det Q ^{\check \delta} (\l)$
is non-zero in a neighbourhood of $\a ^* +i \overline {\a ^* _+}$,
and considering only one Weyl chamber, using Weyl group invariance, we conclude the result.
Using the matrix valued classical Fourier transform (see below)
we conclude that $\H (\goth a ^* _\C) _e ^\delta$ is dense in $\Cal S (\goth a ^* _\v) _e ^\delta$,
and hence $\H (\goth a ^* _\C) _W ^\delta$
is dense in $\Cal S (\goth a ^* _\v) _W ^\delta$.
Elaborating on these results, as in Section 2, we get iii).
$\hfill\square$
\enddemo
Let $\Cal R$ be the Radon transform (replace $a_t \in A$ with $\exp (H),\,H \in \a$).
Let $\p \in C_c ^\infty (\a \times B)$, then the "classical" Fourier transform on $\a \times B$ is defined as:
$$
\Cal F \p (\l ,b) = \int _{\a} \p (H,b) e^{-i\l (H)} dH ,\, \l \in \a ^* _\C ,\, b \in B.
$$
We have the following commutative diagram:
$$
\diagram
& { {\Cal H} ({\a ^* _\C} \times B) ^W}   \\
 { C_c^\infty ( G /K)} \urto^{\Cal H} \rrto_{\Cal R}& &
{ \Cal R C_c^\infty ( G /K) \subset C_c ^\infty (\a \times B)}  \ulto_{\Cal F} \\
\enddiagram
$$As in Section 2, all transforms preserves $K$-types. We define the $\delta$-spherical transform:
$$
\H ^\delta f (\l) = d(\delta) \int _X f(x) \P _{\bar \l , \delta} (x) ^* dx
$$where $\P _{\l, \delta}$ is the generalized sperical function:
$$
\P _{\l,\delta} = \int _K e ^{(i\l +\rho) (A(x,kM))} \delta (k) dk,\, x \in X,
$$
see [He2,\,Chapter\,3,\,\S2,\,\S5].
\proclaim {Theorem 3.3}
The $\delta$-spherical transform $f \mapsto  \Cal H ^\delta f$ is a bijection of $C_c^\infty ( G /K) _{\check
\delta}$ onto $\Cal H (\a ^* _\C) ^\delta _W$.
\endproclaim
\demo {Proof}
[He2,\,Chapter\,III,\,Theorem\,5.11].
$\hfill \square$
\enddemo
Defining an operator $\Cal T$ as in Section 2,
and considering the "classical" Fourier transform on matrix coefficients
$$
\Cal F \p (\l ,b) = \int _{\a} \p (H) e^{i\l (H)} dH ,\, \l \in \a _\C ,
$$
we get the commuting diagram:
$$
\diagram
& {\psi \in {\Cal H} (\a ^* _\C)_W ^\delta }   \\
 {f  \in C_c^\infty ( G /K) _{\check \delta}} \urto^{\Cal H ^\delta}_\cong \rrto_{\Cal T}^\cong& &
{\p \in C_c^\infty (\a) ^\delta _W}  \ulto_{\Cal F}^\cong \\
\enddiagram \tag 4
$$
The subscript $W$ on $C_c^\infty (\a) ^\delta _W$ indicates some kind of symmetry condition.
To show Theorem 3.1 ii)+iii), we are
left to show that ${(\H ^\delta)}^{-1}$ is a continuous homomorphism from
${\Cal H} (\a _\C) _W ^\delta$ into $C_c^\infty ( G /K) _{\check \delta}$ in the induced topologies.
Again we want to use a cut-off function, as in Section 2 and [An]. To do
so, we need to reformulate the Paley-Wiener Theorems for the various transforms.
Define convex W-invariant sets in $\a$ and $G$ by: $\a _R = \{ H \in \a |\r (w\c H ) \le R ,\,\forall w \in W \}$,
and $G_R = K (\exp \a _R )K$. Furthermore define the gauge associated
to $\a _R$ by $ q (\l) = \sup _{H \in \a _R} \l (H)$.
\proclaim {Theorem 3.4}
The "classical" Fourier transform $\F$ is an isomorphism between the space of all functions
$\p \in C _c ^\infty (\a)$, such that $\supp\, \p\subset \a _R$, and the space of all entire functions
$\psi$ on $\a ^* _\C$, such that
$$
\sup _{\l \in \a ^* _\C} e ^{-  q (Im \l ) }(1+\|\l \| ) ^N |\psi (\l)| < \infty, \tag 5
$$
for all $N \in \N$.
\endproclaim
\proclaim {Theorem 3.5} The $\delta$-spherical transform $\H ^\delta$ is an isomorphism between the space of functions $f \in C ^\infty _c (G/K)$, with $\supp f \subset G_R$, and the space of functions $\psi \in \H (\goth a ^* _\C) _W ^\delta $, with matrix entries satisfying $(5)$.
\endproclaim
\demo {Proof}
A slight modification of the proof of [He2,\,Chapter\,III,\,Theorem 5.11], using the ideas from [An,\,p.\,341-2].
$\hfill\square$
\enddemo
\proclaim {Corollary 3.6}
The transform $\Cal T$ is an isomorphism between the spaces
$C ^\infty _c (G/K)$ and $C _c ^\infty (\a) ^\delta _W$.
Moreover $\supp f \subset G_R$ if and only if $\supp \Cal T f \subset \goth a _R$.
\endproclaim
Consider the commuting diagram $(4)$. Given a seminorm $\s$ in $\Cal S ^p (G/K)$,
we shall find a seminorm $\tau$ on $\Cal S (\goth a ^* _\v) _W ^\delta$, such that $\s (f) \le \tau (\psi)$,
for all $f$ and $\psi$.
Instead of looking at the intervals $[0,j]$ and $]j,j+1]$, as we did in Section 2, we will consider
the sets $G_j$ and $G_{j+1}\backslash G_j$. The crucial point in the proof is then to estimate
$f(D;g;E)$ on $G_{j+1}\backslash G_j$.
Let $\omega \in C  ^\infty (\R)$, with $ \omega \equiv 0 $ on $]-\infty ; 0]$,
and $ \omega \equiv 1 $ on $[1 ; \infty[$.
Introduce a $W$-invariant function in $C _c ^\infty (\a)$ by:
$$
\omega _j (H) = \prod _{w \in W} \omega (j - \r (w \c H)).
$$
We see that $\omega _j \equiv 1$ on $\a _{j-1}$,
and  $\omega _j \equiv 0$ outside $\a _{j}$.
Moreover $ \omega _j$ and all its derivatives are bounded in $j$. We decompose $\p$ as:
$$
\split
\p (H)= &\,Q^{\check\delta} \left ( -i\frac {\partial}{\partial H} \right ) \{ (1 - \omega _j ) (H)
\Cal F ^{-1}\{(Q^{\check\delta}) ^{-1} (\c)  \psi (\cdot ) \} (H)\}\\
&+ Q^{\check\delta} \left ( -i\frac {\partial}{\partial H } \right ) \{  \omega _j  (H) \Cal F ^{-1}\{ (Q^{\check\delta}) ^{-1} (\c) \psi (\cdot )  \} (H)\}. \endsplit\
$$
Consider the functions
$$
\split
\p _j(H) &=Q^{\check\delta} \left ( -i\frac {\partial}{\partial H } \right ) \{ (1 - \omega _j ) (H)
\Cal F ^{-1}\{(Q^{\check\delta})^{-1} (\c)  \psi (\cdot ) \} (H)\},\\
\psi _j (\l)  &= {\Cal F \p _j} (\l) = Q^{\check\delta} (\l) \{ \Cal F \{ (1 - \omega _j ) (\cdot) \Cal F ^{-1} \{
(Q^{\check\delta} ) ^{-1} (\c) \psi (\cdot )  \} \} \} (\l), \\
f _j (x) &= ((\Cal H ^\delta ) ^{-1} \psi _j )(x) = (\Cal T ^{-1} \p _j) (x).\\
\endsplit\
$$
We see that $\psi _j \in \H (\goth a ^* _\C) _W ^\delta$ and $\p _j \in C_c^\infty (\a) ^\delta _W$.
Since $\omega _j$ has support in $\a _j$, Corollary 3.6 tells us that $f$
may differ from $f_j$ only inside $G_j$. Continuing as in [An] or Section 2,
using elementary Fourier analysis, we see that the
remaining problem is to estimate $(Q^{\check\delta}) ^{-1} (\c) \psi (\cdot )$ and its derivatives
on the boundary of $\a ^* _\v$, with similar estimates on $\psi (\cdot )$. Following the
argument of Theorem 2.3.9, using the knowledge from the proof of Theorem 3.2, we reach the result.


\Refs
\ref \key {An} \by J.-Ph. Anker \paper The Spherical Fourier Transform of Rapidly Decreasing Functions. A Simple Proof of a Characterization due to Harish-Chandra, Helgason, Trombi and Varadarajan \jour J. Funct. Anal. \vol 96 \yr 1991 \pages 331--349
\endref
\ref \key {Eg}\by M. Eguchi \paper Asymptotic Expansions of Eisenstein Integrals and Fourier Transforms of Symmetric Spaces \jour J. Funct.
Anal.
\vol 34 \yr 1979 \pages 167--216
\endref
\ref \key {GV} \by R. Gangolli and V. S. Varadajaran \paper Harmonic Analysis of Spherical Functions on Real Reductive Groups \publ Ergebnisse der Mathematik
und ihrer Grenzgebiete. Vol. 101, Springer Verlag, Berlin/New York, 1988
\endref
\ref \key {He1} \by S. Helgason \paper Groups and Geometric Analysis \publ Academic Press, Orlando, 1984
\endref
\ref \key {He2} \by S. Helgason \paper Geometric Analysis on Symmetric Spaces \publ Mathematical Surveys and Monographs. Vol. 39, American Mathematical
Society, Providence, Rhode Island, 1994
\endref
\ref \key {Ru} \by W. Rudin \paper Functional Analysis \publ McGraw-Hill, New York, 1973
\endref
\endRefs
\enddocument